\documentclass{IEEEtran}

\usepackage{graphicx}
\usepackage{amssymb}
\usepackage{amsmath}
\usepackage{subfigure}
\usepackage{epsfig}
\usepackage{color}
\usepackage{enumitem}
\usepackage{float}
\usepackage{soul}
\usepackage{wrapfig}
\usepackage{arydshln}

\usepackage{amsthm}
\usepackage [latin1]{inputenc}
\usepackage{comment}

\def\0{{\bf 0}}
\def\1{{\bf 1}}

\def\beq{\begin{equation*}}
    \def\eeq{\end{equation*}}
\def\bql{\begin{equation}}
    \def\eql{\end{equation}}
\def\bqn{\begin{eqnarray*}}
    \def\eqn{\end{eqnarray*}}
\def\bnl{\begin{eqnarray}}
    \def\enl{\end{eqnarray}}
\def\bma{\begin{bmatrix}}
    \def\ema{\end{bmatrix}}
\def\bmx{\begin{matrix}}
    \def\emx{\end{matrix}}
\def\ben{\begin{enumerate}}
    \def\een{\end{enumerate}}
\def\bit{\begin{itemize}}
    \def\eit{\end{itemize}}
\def\bei{\begin{itemize}}
    \def\eei{\end{itemize}}
\def\bet{\begin{tabular}}
    \def\eet{\end{tabular}}

\newcommand{\ba}{\mathbf{a}}

\newcommand{\R}{\mathbb{R}}

\newcommand{\A}{\mathcal{A}}

\newcommand{\be}{\mathbf{e}}

\newcommand{\g}{\mathbf{g}}

\newcommand{\bu}{\mathbf{u}}

\newcommand{\bv}{\mathbf{v}}

\usepackage[dvipsnames]{xcolor}

\def\R{\mathbb{R}}

\def\1{{\bf1}}
\def\la{\langle}
\def\ra{\rangle}
\def\b{{\beta}}
\def\a{\alpha}
\def\g{\gamma}

\def\bit{\begin{itemize}}
\def\eit{\end{itemize}}

\def\be{\begin{equation}}
\def\ee{\end{equation}}
\def\ba{\begin{eqnarray}}
\def\ea{\end{eqnarray}}

\def\bes{\begin{equation*}}
\def\ees{\end{equation*}}
\def\bas{\begin{eqnarray*}}
\def\eas{\end{eqnarray*}}

\usepackage{url}

\newtheorem{Remark 1}{Remark}
\newtheorem{Remark 2}[Remark 1]{Remark}
\newtheorem{Remark 3}[Remark 1]{Remark}
\newtheorem{Remark 4}[Remark 1]{Remark}
\newtheorem{Remark 5}[Remark 1]{Remark}
\newtheorem{Remark 6}[Remark 1]{Remark}
\newtheorem{Remark 7}[Remark 1]{Remark}

\newtheorem{Lemma 1}{Lemma}
\newtheorem{Lemma 2}[Lemma 1]{Lemma}
\newtheorem{Lemma 3}[Lemma 1]{Lemma}
\newtheorem{Lemma 4}[Lemma 1]{Lemma}
\newtheorem{Lemma 5}[Lemma 1]{Lemma}
\newtheorem{Lemma 6}[Lemma 1]{Lemma}
\newtheorem{Lemma 7}[Lemma 1]{Lemma}

\newtheorem{Proposition 1}{Proposition}
\newtheorem{Proposition 2}[Proposition 1]{Proposition}

\newtheorem{Assumption 1}{Assumption}
\newtheorem{Assumption 2}[Assumption 1]{Assumption}
\newtheorem{Assumption 3}[Assumption 1]{Assumption}
\newtheorem{Assumption 4}[Assumption 1]{Assumption}
\newtheorem{Definition 1}{Definition}
\newtheorem{Theorem 1}{Theorem}
\newtheorem{Theorem 2}[Theorem 1]{Theorem}
\newtheorem{Theorem 3}[Theorem 1]{Theorem}
\newtheorem{Theorem 4}[Theorem 1]{Theorem}
\newtheorem{Theorem 5}[Theorem 1]{Theorem}
\newtheorem{Theorem 6}[Theorem 1]{Theorem}
\newtheorem{Theorem 7}[Theorem 1]{Theorem}
\newtheorem{Theorem 8}[Theorem 1]{Theorem}
\newtheorem{Theorem 9}[Theorem 1]{Theorem}
\newtheorem{Theorem 10}[Theorem 1]{Theorem}

\title{\LARGE \bf
 Decentralized Gradient Methods with Time-varying Uncoordinated Stepsizes: Convergence Analysis and Privacy Design}

\author{Yongqiang Wang, Angelia Nedi\'c
\thanks{ The work was supported in part by the National Science Foundation under Grants ECCS-1912702, CCF-2106293, CCF-2106336, CCF-2215088,
and CNS-2219487.}
\thanks{Yongqiang Wang is with the Department of Electrical and Computer Engineering, Clemson University, Clemson, SC 29634, USA
{\tt\small{yongqiw}@clemson.edu}
}%
\thanks{Angelia Nedi\'c is with the School of Electrical, Computer and Energy Engineering, Arizona State University, Tempe, AZ 85281, USA {\tt\small angelia.nedich@asu.edu}}
}

\begin{document}

\maketitle
\thispagestyle{empty}
\pagestyle{empty}

\begin{abstract}
Decentralized optimization enables a network of agents to
cooperatively optimize an overall objective function without a
central coordinator and is gaining increased attention in domains as
diverse as control, sensor networks, data mining, and robotics.
However, the information sharing among agents in decentralized
optimization also  discloses  agents' information, which is
undesirable or even unacceptable when involved data are  sensitive.
This paper proposes two gradient based decentralized optimization
algorithms that can protect participating agents' privacy without
compromising optimization accuracy or incurring heavy
communication/computational overhead. This is in distinct difference
from differential privacy based approaches which have to trade
optimization accuracy  for privacy, or encryption based approaches
which incur heavy communication and computational overhead. Both
algorithms leverage a judiciously designed mixing matrix and
time-varying uncoordinated stepsizes to enable privacy, one using
diminishing stepsizes while the other   using non-diminishing
stepsizes.  In both algorithms, when interacting with any one of its neighbors,   a participating agent  only needs to
share {\it one} message  in each iteration
to reach convergence  to an exact optimal solution, which is in
contrast to most gradient-tracking based algorithms requiring every
agent to share two messages (an optimization variable and an
auxiliary gradient-tracking variable)
  under non-diminishing stepsizes. Furthermore, both algorithms
can guarantee the privacy of a participating agent even when all
information  shared by the agent are accessible to an adversary, a
scenario in which
  most existing accuracy-maintaining privacy  approaches will fail to protect privacy. Simulation results confirm the effectiveness of the proposed
algorithms.
\end{abstract}

\section{Introduction}
Distributed optimization is gaining increased attention across
disciplines due to its fundamental importance and vast applications
in areas ranging from cooperative control  \cite{yang2019survey},
distributed sensing \cite{bazerque2009distributed}, multi-agent
systems \cite{raffard2004distributed}, sensor networks
\cite{zhang2017distributed}, to large-scale machine learning
\cite{tsianos2012consensus}. In many of these applications, the
problem can be formulated in the following general form, in which a
network of $m$ agents cooperatively solve a common optimization
problem through on-node computation and local communication:
\begin{equation}\label{eq:optimization_formulation1}
\min\limits_{\theta\in\mathbb{R}^d} F(\theta)\triangleq
\frac{1}{m}\sum_{i=1}^m f_i(\theta)
\end{equation}
where $\theta$ is  common to all agents but
$f_i:\mathbb{R}^d\rightarrow\mathbb{R}$ is a local objective
function private to agent $i$. We denote an optimal solution to this
problem by $\theta^{\ast}$, which we assume to be finite.

Since  the 1980s,
the above decentralized optimization problem has been intensively
studied. To date, various algorithms have been proposed. Some of
the  commonly used algorithms to solve
(\ref{eq:optimization_formulation1})   include decentralized gradient methods (e.g.,
\cite{nedic2009distributed,shi2015extra,xu2015augmented,xin2018linear}),
distributed alternating direction method of multipliers (e.g.,
\cite{shi2014linear,zhang2019admm}),  and distributed Newton methods
(e.g., \cite{wei2013distributed}). We focus on the gradient based
approach due to its simplicity in computation, which is particularly
appealing when   agents have limited computational
capabilities.

Over the past decade, plenty of gradient based  algorithms have been
developed for decentralized optimization. Early results combine
consensus  and gradient method  by directly concatenating gradient
based step  with a  consensus operation of the optimization
variable. Typical examples include
\cite{nedic2009distributed,yuan2016convergence}. However, to find an
exact optimal solution, these approaches have to use a diminishing
stepsize, which slows down the convergence.
To guarantee both a fast
convergence speed and exact optimization result, algorithms have
been proposed to replace the local gradient in decentralized
gradient methods with an auxiliary variable which tracks the
gradient of the global objective function. Typical examples include
Aug-DGM \cite{xu2015augmented}, DIGing \cite{nedic2017achieving},
AsynDGM \cite{xu2017convergence}, AB \cite{xin2018linear}, Push-Pull
\cite{pu2020push,du2018accelerated} and ADD-OPT \cite{xi2017add},
etc. While these algorithms can converge to an exact optimal
solution under a fixed stepsize,   they   have to exchange both the optimization
variable and the additional auxiliary variable
 in every iteration, which
{\it doubles} the communication overhead in each iteration compared
with  conventional decentralized gradient based algorithms.

All of the aforementioned algorithms explicitly share optimization
variables and/or gradients in every iteration, which  becomes a
problem in applications involving sensitive data. For example, in
the rendezvous problem where a group of agents use decentralized
optimization to cooperatively find an optimal assembly point,
participating agents may want to keep their initial positions
private, which is particularly important in unfriendly environments
\cite{zhang2019admm}. In fact, without an effective privacy
mechanism in place, the results in
\cite{zhang2019admm,huang2015differentially,burbano2019inferring}
show  that a participating agent's sensitive information such as
position can be easily  inferred by an adversary or other
participating agents in decentralized-optimization based rendezvous
and parameter estimation. Another case underscoring the importance
of privacy preservation in decentralized optimization is machine
learning where exchanged data may contain sensitive information such
as medical records and salary information \cite{yan2012distributed,wei2020federated}.
In fact, recent  results in \cite{zhu2019deep} show  that without a
privacy mechanism, an adversary can precisely recover the raw data
(pixel-wise accurate for images and token-wise matching for texts)
through shared gradients.

Recently results have been reported on
privacy-preserving decentralized optimization. For example,
differential-privacy based approaches have been proposed to obscure
shared information in decentralized optimization by injecting noise
to exchanged messages
\cite{huang2015differentially,cortes2016differential,xiong2020privacy,wang2022quantization,wang2022tailoring}
or objective functions \cite{nozari2016differentially}. However, the added noise in differential
privacy also unavoidably compromises the accuracy of optimization
results. To enable privacy protection without sacrificing
optimization accuracy, partially homomorphic encryption has been
employed in both our own prior results
\cite{zhang2019admm,zhang2018enabling}, and others'
\cite{freris2016distributed,lu2018privacy}. However, such approaches
   incur  heavy communication and computation overhead. Employing the structural
  properties of decentralized optimization, results have also been
  reported on privacy protection in decentralized optimization without
  using differential privacy or encryption. For example,
  \cite{yan2012distributed,lou2017privacy} showed that privacy can be enabled by
  adding a {\it constant} uncertain parameter in the projection step or stepsizes. The authors of \cite{gade2018private} showed that network structure can be leveraged to
  construct spatially correlated ``structured" noise to cover
  information.
   Although these
  approaches can ensure convergence to an exact  optimal solution,
  their enabled privacy is restricted: projection based privacy depends on
  the size of the projection set -- a large projection set nullifies
  privacy protection whereas a small projection set offers strong
  privacy
  protection but requires {\it a priori} knowledge of the optimal solution; ``structured"
  noise based approach requires each agent to have a certain number
  of neighbors which do not share information with the adversary.  In
  fact, such a structure constraint is required in most existing
  privacy solutions with guaranteed optimization accuracy. For example,
   besides the ``structured" noise based approach, even the partially homomorphic encryption based privacy
  approaches for decentralized optimization require the adversary not  to have access to a target agent's communications with all of its neighbors  \cite{zhang2019admm}.

  Inspired by our recent results  that privacy can be
  enabled in consensus by manipulating inherent dynamics
  \cite{ruan2019secure,wang2019privacy,gao2022algorithm},  we propose to
 enable privacy in decentralized gradient methods  by judiciously manipulating the inherent
 dynamics of information mixing and gradient  operations.
 More specifically, leveraging a judiciously designed mixing matrix and time-varying uncoordinated stepsizes, we propose two privacy-preserving decentralized gradient based algorithms, one
 with
 diminishing stepsizes and the other one with non-diminishing stepsizes.
 Not only do  our algorithms maintain the accuracy of
 decentralized optimization, they also
  enable privacy even when an adversary has access to all
  messages shared by a participating agent. This is in contrast to most
  existing accuracy-guaranteed privacy approaches for decentralized
  optimization which cannot protect an agent against adversaries having
  access to all shared messages.
   Furthermore, even in the non-diminishing stepsize case, our algorithm only requires a participating agent to
   share  one variable with { any one of its neighboring agents} in each iteration, which is extremely appealing when
   communication bandwidth is limited. In fact, to our
   knowledge, our algorithm is the first   privacy-preserving
   decentralized gradient based algorithm  that uses non-diminishing stepsizes to reach accurate optimization results but requires each participating  agent to
   share only  one message with a neighboring agent
   in every iteration. Note that most existing
   gradient-tracking based decentralized optimization algorithms (e.g., \cite{xu2015augmented,xin2018linear,xu2017convergence,pu2020push,du2018accelerated,xi2017add,di2016next,daneshmand2020second,xin2020general,bin2019system,zhang2020distributed}) require an agent to share
   two messages (the optimization variable and an auxiliary variable tracking the gradient of the global objective function)  in every iteration.

   The main contributions are as follows: 1)~We
   propose two accuracy-guaranteed decentralized gradient based
   algorithms that can protect the privacy of  participating agents even when all  shared messages are accessible to
   an adversary, a scenario which fails   existing
   accuracy-guaranteed privacy-preserving approaches for decentralized
   optimization. {The first one adapts a privacy-preserving algorithm for   decentralized stochastic optimization in \cite{wang2022decentralized_inherent} and the second one is a brand new algorithm}; 2) The two inherently privacy-preserving algorithms
   are efficient in communication/computation  in that they are
   encryption-free and only require a participating agent to share {\it one} message with a neighboring agent in every iteration, both   for  the diminishing stepsize case (the first algorithm)
    and
    the non-diminishing stepsize case (the second algorithm). This is significant in that, as a comparison,   existing
   state-of-the-art gradient-tracking based decentralized optimization
   algorithms require a participating agent to share two messages in every iteration, both the optimization variable and the auxiliary gradient-tracking variable\footnote{Some
    gradient-tracking based decentralized optimization algorithms
   may be transformed to and implemented in a one-variable form --like EXTRA \cite{shi2015extra}-- and share one message in every iteration; however, such an implementation becomes
   infeasible when the stepsizes or coupling weights are
   uncoordinated to enable privacy. See  Remark \ref{remark:conversion} and Sec. \ref{sec:privacy_analysis_NDS} for detailed explanations.}. In fact, the
   sharing of the additional gradient-tracking variable will lead to serious privacy
   breaches, as  detailed in Sec. \ref{sec:privacy_analysis_NDS}; 3) Even without taking privacy into
   consideration, the two decentralized optimization algorithms
   are of interest by themselves. To our knowledge, our
   analysis is the first to rigorously characterize the
   convergence of decentralized gradient methods in the presence of
   time-varying heterogeneity in  stepsizes, which is in contrast to existing results only
   addressing
   constant or fixed
 heterogeneity in stepsizes
 \cite{xu2015augmented,lou2017privacy,nedic2017geometrically}.

The organization of the paper is as follows. Sec. II gives the
problem formulation. Sec. III presents PDG-DS, an inherently
privacy-preserving decentralized gradient algorithm with proven
converge to the accurate optimization solution under diminishing
uncoordinated  stepsizes. Both its convergence properties and
privacy performance are also rigorously characterized. Sec. IV
presents PDG-NDS, an inherently privacy-preserving decentralized
gradient based algorithm with proven converge  to the exact
optimization solution under non-diminishing and time-varying
uncoordinated stepsizes. Its convergence properties and privacy
performance are also rigorously characterized in this section.
Subsection \ref{sec:privacy_analysis_NDS} of this section also
analyzes why existing gradient-tracking based decentralized
optimization algorithms   breach the privacy of participating
agents, even when coupling weights are time-varying and
 stepsizes are uncoordinated. Sec. V gives simulation results as well as
comparison with existing works. Finally Sec. VI concludes the paper.

{\bf Notations:} $\mathbb{R}^m$ denotes the Euclidean space of
dimension $m$. $I_d$ denotes identity matrix of dimension $d$. ${\bf
1}_d$ denotes a $d$ dimensional column vector will all entries equal
to 1 and we omit the dimension when clear from the context. A vector
is viewed as a column vector. For a vector $x$, $x_i$ denotes the
$i$th element. We say $x> 0$ (resp. $x\geq 0$) if all elements of
$x$ are positive (resp. non-negative). $A^T$ denotes the transpose
of matrix $A$ and $\langle x, y\rangle$ denotes the inner product of
two vectors $x$ and $y$. $\|\cdot\|$ denotes the standard Euclidean
norm  { for a vector or the induced Euclidean norm (spectral  norm) for a matrix}. A  matrix  is
column-stochastic when its  entries are nonnegative and  elements in
every
column add up to one. { A matrix is doubly-stochastic when its entries are nonnegative and  elements in
every
column   add up to one, and elements in every row add up to one.}  For two matrices $A$ and $B$ with the same dimensions, we use $A\leq B$ to represent that every entry of $A$ is no larger than the corresponding entry of $B$.   We use $\otimes$ to represent the Kronecker product. 

\section{Problem Formulation}

We consider a network of $m$ agents. The agents interact on an
undirected graph, which can be described by a weight matrix
$W=\{w_{ij}\}$. More specifically, if agents  $i$ and  $j$ can
  interact with each other, then $w_{ij}$ is
positive. Otherwise, $w_{ij}$ will be zero. We assume that an agent
is always able to affect itself, i.e., $w_{ii}>0$ for all $1\leq
i\leq m$. The neighbor set $\mathbb{N}_i$ of agent $i$ is defined as
the set of agents $\{j|w_{ij}>0\}$. So the neighbor set of agent $i$
always includes itself.

\begin{Assumption 1}\label{assumption:W}
   $W=\{w_{ij}\}\in \mathbb{R}^{m\times m}$
     satisfies ${\bf 1}^TW={\bf
  1}^T$, $W{\bf 1}={\bf
  1}$, and $\eta=\|W-\frac{{\bf 1}{\bf 1}^T}{m}\|<1$.
\end{Assumption 1}

The optimization problem (\ref{eq:optimization_formulation1}) can be
reformulated as the following equivalent multi-agent optimization
problem:
\begin{equation}\label{eq:optimization_formulation2}
\min\limits_{x\in\mathbb{R}^{md}}f(x)\triangleq
\frac{1}{m}\sum_{i=1}^m f_i(x_i)\: {\rm s.t.}\:
x_1=x_2=\cdots=x_m
\end{equation}
where $x_i\in\mathbb{R}^d$ is the local estimate of agent $i$ about
the optimization solution and
$x=[x_1^T,x_2^T,\cdots,x_m^T]^T\in\mathbb{R}^{md}$ is the collection
of the estimates of all agents.

 We make the following standard
assumption on  objective functions:
\begin{Assumption 2}\label{assumption:L_and_G}
Problem (\ref{eq:optimization_formulation1}) has at least one
optimal solution $\theta^{\ast}$. Every $f_i$  has
Lipschitz continuous gradients, i.e., for some
$L>0$, $\|\nabla f_i(u)-\nabla f_i(v)\|\le L \|u-v\|,\quad\forall i \:\:{\rm and}\:\:\forall u,v\in\mathbb{R}^d.$ Every $f_i$ is {  convex, i.e., $f_i(u)\geq f_i(v)+\nabla f_i(v)^T(u-v) \: \forall i \: {\rm and}\:\:\forall u,v\in\mathbb{R}^d$}.
\end{Assumption 2}

Under Assumption \ref{assumption:L_and_G}, we know that
(\ref{eq:optimization_formulation2})   always has an optimal
solution $ x^{\ast}=
[(\theta^{\ast})^T,(\theta^{\ast})^T,\cdots,(\theta^{\ast})^T]^T $.

In decentralized optimization applications, gradients usually carry
sensitive information. For example, in decentralized-optimization
based rendezvous and localization, disclosing the gradient  of an
agent amounts to disclosing its (initial) position
\cite{zhang2019admm,huang2015differentially}. In machine learning,
gradients are directly calculated  from and embed information of
sensitive learning data \cite{zhu2019deep}. Therefore, in this
paper, we define privacy as preventing disclosing agents' gradients { in each  iteration.}

  We consider two potential attacks, which are the
two most commonly used models of attacks in privacy research
\cite{Goldreich_2}:
\begin{itemize}
\item \emph{Honest-but-curious attacks}  are attacks in which a participating
agent or multiple participating agents (colluding or not)
 follows all protocol steps correctly but is curious and
collects all received intermediate data  to learn the sensitive
information about other participating agents.

\item \emph{Eavesdropping attacks}  are attacks in which an external eavesdropper wiretaps all communication channels to
intercept exchanged messages so as to learn sensitive information
about sending agents.
\end{itemize}

{  An honest-but-curious adversary (e.g., agent $i$) has access
to the internal state $x_i$, which is unavailable to external eavesdroppers. However, an eavesdropper has access to all shared
information in the network, whereas an honest-but-curious agent
can only access shared information that is destined to it.}

\section{An inherently privacy-preserving decentralized gradient algorithm with diminishing stepsizes}

Conventional decentralized gradient algorithms usually take the
following form:
\begin{equation}\label{eq:conventional_gradient_descent}
 x_i^{k+1}=\sum\nolimits_{j\in\mathbb{N}_i}w_{ij}x_j^k-\lambda^k g_i^k
\end{equation}
where $\lambda^k$ is a positive scalar denoting the stepsize and
  $g^k_i$ denotes
the gradient of agent $i$ evaluated at $x_i^k$. It is well-known
that under Assumption \ref{assumption:W} and Assumption
\ref{assumption:L_and_G}, when $\lambda^k$ diminishes with $k$,
i.e., $\sum_{k=0}^{\infty}\lambda^k=\infty$ and
$\sum_{k=0}^{\infty}(\lambda^k)^2<\infty$, then all $x_i^k$ will
converge to a same optimal solution. {  It is worth noting that by diminishing stepsizes, we mean  stepsizes that eventually converge to zero. This does not require the stepsize to decrease monotonically. }

However, in the above decentralized gradient  algorithm, agent $i$
has to share $x^k_i$
 with all its neighbors $j\in\mathbb{N}_i$.  If
an adversary has access to the optimization variable $x_i^k$ of
agent $i$
 and the updates that agent $i$ receives from all its neighbors
$x_j^k$ for $j\in\mathbb{N}_i$, then the adversary can easily infer
$g_i^k$   based on the update rule
(\ref{eq:conventional_gradient_descent}) and publicly known  $W$ and
$\lambda^k$.

Motivated by this observation and inspired by our recent finding
that interaction dynamics can be judiciously manipulated to enable
privacy \cite{ruan2019secure,wang2019privacy,gao2022algorithm}, we propose the
following  decentralized gradient  algorithm to enable privacy (with
per-agent version given in Algorithm PDG-DS on the next page) { by adapting the decentralized stochastic optimization algorithm in \cite{wang2022decentralized_inherent}}:
\begin{equation}\label{eq:decaying_stepsize}
x^{k+1}=(W\otimes I_d)x^k-((B^k\Lambda^k)\otimes I_d)g^k
\end{equation}
where $B^k=\{b_{ij}^k\}$ is a column-stochastic nonnegative matrix,
$\Lambda^k={\rm diag}[\lambda_1^k,\lambda_2^k,\cdots,\lambda_m^k]$
with { $\lambda_i^k\geq 0$} denoting the stepsize of agent $i$ at
iteration $k$ and $g^k=[(g_1^k)^T,(g_2^k)^T,\cdots (g_m^k)^T]^T$. It is worth noting that different from the algorithm in \cite{wang2022decentralized_inherent} which uses a matrix-valued stepsize (with dimension $\mathbb{R}^{d\times d}$) for each agent, we here require the stepsize $\lambda_j^k$ for each agent to be a scalar, which is necessary in our derivation to prove deterministic convergence to an exact optimal solution. { It is also worth noting that our scalar stepsize  here cannot be viewed as a special case of the matrix-valued stepsize in \cite{wang2022decentralized_inherent}  since  the stepsize matrix in \cite{wang2022decentralized_inherent} explicitly requires all diagonal entries  to be  statistically independent of each other, which  prohibits it from having the form of $\lambda_j^k I_d$ (in which all diagonal entries are equal to $\lambda_j^k$ and hence are statistically  dependent on each other).}

 The detailed implementation procedure for individual agents is
provided in Algorithm PDG-DS. Compared with the conventional
decentralized gradient  algorithm, it can be seen that we equip
each agent $j$ with  two private variables  $b^k_{ij}$ and
$\lambda_j^k$ to cover its gradient information $g_j^k$. The two
variables are generated by and only known to agent $j$. Therefore,
the two variables can ensure that a neighboring agent $i$ cannot
infer $g_j^k$ based on received information
$w_{ij}x_j^k-b_{ij}^k\lambda_j^kg_j^k$ as agent $i$ does not know
$x_j^k$, $\lambda_j^k$, or $b_{ij}^k$. { It is worth noting that since
  agent $j$ determines $b_{ij}^k>0$ for all $i\in\mathbb{N}_j$ ($b_{ij}^k=0$ for  $i\notin\mathbb{N}_j$), as
long as every agent $j$ ensures $\sum_{i\in\mathbb{N}_j}b_{ij}^k=1$
locally, the column-stochastic condition for matrix $B^k$ will be
satisfied.} Furthermore, as every agent $j$ determines its stepsize
$\lambda_j^k$ randomly and independently from each other at each
iteration, the stepsize will be heterogeneous across the agents and
the difference between stepsizes of two agents would also be varying
with time.

In the following, we will first prove that even under the
time-varying mixing matrix $B^k$ and time-varying heterogeneous
stepsizes, Algorithm PDG-DS can still guarantee convergence of all
agents to an exact optimal solution. Then we will rigorously analyze
its performance of privacy protection. It is worth noting that
adding $B^k$ will not only facilitate privacy protection, it will
also enhance information mixture across agents and hence speed up
convergence to the optimal solution, as illustrated in the numerical
simulations in Sec. V.

 \noindent\rule{0.49\textwidth}{0.5pt}
\noindent\textbf{PDG-DS: Privacy-preserving decentralized gradient
method with diminishing stepsizes}

\vspace{-0.2cm}\noindent\rule{0.49\textwidth}{0.5pt}
\begin{enumerate}
    \item[] Public parameters: $W$
    \item[] Private parameters for each agent $i$: $b^k_{ji}$, $\lambda_i^k$, and {$x_i^0$}
    \item {\bf for  $k=0,1,\cdots$ do}
    \begin{enumerate}
        \item Every agent $j$ computes and sends to  agent
        $i\in\mathbb{N}_j$
                \begin{equation}\label{eq:v_ij_DS}
             v^k_{ij}\triangleq { w_{ij}}x_j^k-b_{ij}^k\lambda_j^kg_j^k
        \end{equation}
        \item After receiving $v_{ij}^k$ from all $j\in\mathbb{N}_i$, agent $i$ updates its state as follows:
        \begin{equation}\label{update-rule-01}
            x_i^{k+1}=\sum\nolimits_{j\in\mathbb{N}_i}v_{ij}^k=\sum\nolimits_{j\in\mathbb{N}_i}({w_{ij}}x_j^k-b_{ij}^k\lambda_j^kg_j^k)
        \end{equation}
        \vspace{-0.4cm}
        \item {\bf end}
    \end{enumerate}
\end{enumerate}
\vspace{-0.2cm}\rule{0.49\textwidth}{0.5pt}
\subsection{Convergence Analysis}
For the convenience of analysis, we   define the average vector
$\bar x^k=\frac{\sum_{i=1}^m x_i^k}{m}$. Because $B^k$ and $W$ are
column stochastic, from (\ref{eq:decaying_stepsize}), we can obtain
\begin{equation}\label{eq:bar_x^k}
\bar x^{k+1}=\bar x^k - \textstyle\frac{1}{m}
\sum_{i=1}^m\lambda_i^k g_i^k
\end{equation}

 To analyze  PDG-DS, we first introduce a theorem that applies to
general decentralized algorithms for solving
(\ref{eq:optimization_formulation1}).
\begin{Proposition 1}\label{th-main}
Assume that problem (\ref{eq:optimization_formulation1}) is convex
and has a solution. Suppose that a distributed algorithm generates
sequences $\{x_i^k\}\subseteq\R^d$ such that the following relation
is satisfied for any optimal solution $\theta^*$ and for all $k$,
\ba\label{eq-fin}  \bv^{k+1} \le \left( \left[\begin{array}{cc} 1 &
1\cr 0& \eta\cr
\end{array}\right]
+a^k {\bf 1}{\bf 1}^T  \right)\bv^k +b^k{\bf 1} - c^k
\left[\begin{array}{c}
 \zeta\cr 0\end{array}\right]
\ea where
$\bv^k=\left[\begin{array}{c}\bv_1^k\\\bv_2^k\end{array}\right]\triangleq
\left[\begin{array}{c}\|\bar x^k-\theta^*\|^2\cr
\sum_{i=1}^m\|x_i^k-\bar x^k\|^2\end{array}\right]$,
 $\zeta\triangleq\sum_{i=1}^m \left(f_i(\bar x^k)- f_i(\theta^*)
\right)$, and the scalar sequences $\{a^k\}$, $\{b^k\}$, and
$\{c^k\}$ are nonnegative and satisfy $\sum_{k=0}^\infty
a^k<\infty$, $\sum_{k=0}^\infty b^k<\infty$, and $\sum_{k=0}^\infty
c^k=\infty$. Then, we have $\lim_{k\to\infty}\|x_i^k - \bar x^k\|=0$
for all  $i$ and there exists an optimal solution
$\tilde\theta^{\ast}$ such that $\lim_{k\to\infty}\|\bar x^k-
\tilde\theta^*\|=0$.
\end{Proposition 1}


\begin{proof}
Since $\theta^*$ is an  optimal solution of problem
(\ref{eq:optimization_formulation1}), we always have $\sum_{i=1}^m
\left(f_i(\bar x^k) -f_i(\theta^*)\right)\ge0$.

From \eqref{eq-fin} it follows that for all $k\ge0$,
\be\label{eq-fin0} \bv^{k+1} \le \left( \left[\begin{array}{cc} 1 &
1\cr 0& \eta\cr\end{array}\right] +a^k {\bf 1}{\bf 1}^T
\right)\bv^k+b^k{\bf 1} \ee Consider the vector $\pi=[1,
\frac{1}{1-\eta}]^T$ and note  $\pi^T \left[\begin{array}{cc} 1 &
1\cr 0& \eta\cr\end{array}\right]= \pi^T$. Thus,
  the sequence $\{\pi^T\bv^k\}$ satisfies all
conditions of Lemma~\ref{lem-polyak} in the Appendix. Therefore, it
follows that
  $\lim_{k\to\infty}\pi^T\bv^k$ exists and that  $\{\|\bar x^k-\theta^*\|\}$ and $\{\sum_{i=1}^m \|x^k-\bar
x^k\|^2\}$ are bounded.

We use  $M>0$ to represent an upper bound on  $\{\|\bar
x^k-\theta^*\|\}$ and $\{\sum_{i=1}^m \|x^k-\bar x^k\|^2\}$, i.e.,
  $\|\bar x^k-\theta^*\|\le
M$ and $\sum_{i=1}^m \|x^k-\bar x^k\|^2\le M$ hold  $\forall k\ge0$.
Thus, for
  all
$k\ge0$, we have
\begin{equation}\label{eq-sumxk}
\begin{aligned}
& \sum_{i=1}^m \|x_i^{k+1}-\bar
x^{k+1}\|^2  \le   \eta\sum_{i=1}^m \|x_i^k - \bar x^k\|^2 + 2 a^k M+b^k\\
& = \sum_{i=1}^m \|x_i^k - \bar x^k\|^2 - (1-\eta)\sum_{i=1}^m
\|x_i^k - \bar x^k\|^2+ 2 a^k M+b^k
\end{aligned}
\end{equation}
By summing (\ref{eq-sumxk}) over $k$ and using the fact
$\sum_{k=0}^\infty (2 a^k M+b^k)<\infty$, we obtain
\be\label{eq-sumable} (1-\eta)\sum\nolimits_{k=0}^\infty
\sum\nolimits_{i=1}^m \|x_i^k - \bar x^k\|^2<\infty\ee which
implies $\lim_{k\to\infty} \|x_i^k - \bar x^k\|^2=0$ for all $i$.

Next, we consider the first element of $\bv^k$, i.e., $\|\bar
x^k-\theta^*\|^2$. From  \eqref{eq-fin} we have
\[
\begin{aligned}\|\bar x^{k+1}-\theta^*\|^2 &\le  (1+a^k)\left(\|\bar x^k
-\theta^*\|^2 +\sum\nolimits_{i=1}^m \|x_i^k - \bar x^k\|^2\right) \\
&\quad + b^k -c^k \sum\nolimits_{i=1}^m\left(f_i(\bar x^k)-
f_i(\theta^*) \right)\cr &\le
 \|\bar x^k -\theta^*\|^2 +\sum\nolimits_{i=1}^m \|x_i^k - \bar x^k\|^2 +2 a^k
M\\
&\quad +b^k - c^k \sum\nolimits_{i=1}^m\left(f_i(\bar x^k)-
f_i(\theta^*) \right)
\end{aligned}
\]

We can see that the preceding relation satisfies the relation in
Lemma~\ref{lem-opt} in the Appendix with $\phi= \sum_{i=1}^m f_i$,
$z^*=\theta^*$,   $z^k =\bar x^k$, $\a^k=0$, $\g^k= c^k$, and
$\beta^k=\sum_{i=1}^m \|x_i^k - \bar x^k\|^2+2 a^k M+b^k$. By our
assumption, we have $\sum_{k=0}^\infty a^k<\infty$,
$\sum_{k=0}^\infty b^k<\infty$, and   $\sum_{k=0}^\infty
c^k=\infty$. Thus, in view of relation~\eqref{eq-sumable}, it
follows $\sum_{k=0}^\infty\beta^k<\infty$. Hence, all   conditions
of Lemma~\ref{lem-opt} in the Appendix are satisfied, and it follows
that $\{\bar x^k\}$ converges to some optimal solution.
\end{proof}

Now, we are in position to prove convergence of PDG-DS.
\begin{Theorem 1}\label{theorem_PDGD_DS_convergence}
Under Assumption \ref{assumption:W} and Assumption
\ref{assumption:L_and_G}, if the stepsize of every agent $i$ {  is non-negative and}
satisfies $\sum_{k=0}^{\infty}\lambda_i^k=+\infty$ and
$\sum_{k=0}^{\infty}(\lambda_i^k)^2<\infty$, and   the stepsize
heterogeneity satisfies
\begin{equation}\label{eq:stepsize_heterogenerity}
 {\sum\nolimits_{k=0}^{\infty}\: \sum\nolimits_{i,j\in\{1,2,\cdots,m\},\,i\neq j}|\lambda_i^{k}-\lambda_j^k|}<\infty
\end{equation}
then we have  $\lim_{k\to\infty}\|x_i^k - \bar x^k\|=0$ for all
  $i$, and there exists an optimal solution
$\tilde\theta^{\ast}$ such that   $\lim_{k\to\infty}\|\bar x^k-
\tilde\theta^*\|=0$ holds.
\end{Theorem 1}
\begin{proof}
The basic idea is to show that  Proposition \ref{th-main} applies. So we
have to establish necessary relationships  for   $\|\bar
x^k-\theta^*\|^2$ and $\sum_{i=1}^m\|x_i^k-\bar x^k\|^2$, which are
fulfilled in Step I and Step II below, respectively.

Step I: Relationship for $\|\bar x^k-\theta^*\|^2$. Using
(\ref{eq:bar_x^k}), we have for any optimal solution $\theta^*$
\[\bar x^{k+1}-\theta^*=\bar x^k -\theta^* - \textstyle\frac{1}{m} {\sum_{i=1}^m\lambda_i^k g_i^k} \]
which further implies
\begin{equation}\label{eq-rel1}
\begin{aligned}
 \left\|\bar x^{k+1}-\theta^*\right\|^2=&\left\|\bar
x^k -\theta^* \right\|^2 - \textstyle\frac{2}{m} \sum_{i=1}^m
\langle \lambda_i^kg_i^k, \bar x^k -\theta^*
\rangle  \\
 &+ \textstyle\frac{1}{m^2} \left\|\sum_{i=1}^m\lambda_i^k
 g_i^k\right\|^2
\end{aligned}
\end{equation}
We next estimate the inner product term, for which we have
\begin{equation}\label{eq:inner_product_0}
\begin{aligned}
 \langle \lambda_i^kg_i^k, \bar x^k -\theta^* \rangle= &\langle \lambda_i^k(g_i^k-\nabla f_i(\bar{x}^k)),  \bar x^k -\theta^*
 \rangle\\
 &+\langle\lambda_i^k\nabla f_i(\bar{x}^k), \bar x^k -\theta^*
 \rangle
 \end{aligned}
\end{equation}
By  the Lipschitz continuous property of $\nabla f_i$, we obtain
\begin{equation}\label{eq:inner_product_1}
\begin{aligned}
&\langle \lambda_i^k(g_i^k-\nabla f_i(\bar{x}^k)),  \bar x^k
-\theta^*
 \rangle\geq -L\lambda_i^k\|x_i^k-\bar{x}^k\|\|\bar x^k -\theta^*\|\\
&\geq-\textstyle\frac{1}{2}\|x_i^k-\bar{x}^k\|^2
-\textstyle\frac{1}{2}L^2(\lambda_i^k)^2\|\bar x^k -\theta^*\|^2
\end{aligned}
\end{equation}
Defining the average stepsize
$\bar{\lambda}^k=\frac{\sum\lambda_i^k}{m}$, we have
\begin{equation}\label{eq:inner_product_2}
\begin{aligned}
\hspace{-0.2cm}&\langle\lambda_i^k\nabla f_i(\bar{x}^k), \bar x^k
-\theta^*
 \rangle=\\
\hspace{-0.2cm} &\langle(\lambda_i^k-\bar\lambda^k)\nabla
f_i(\bar{x}^k), \bar x^k -\theta^*
 \rangle +\langle\bar\lambda^k\nabla f_i(\bar{x}^k), \bar x^k -\theta^*
 \rangle
 \end{aligned}
\end{equation}

Defining $\lambda^k = [\lambda_1^k, \cdots,\lambda_m^k]^T$ and
combining (\ref{eq:inner_product_0})-(\ref{eq:inner_product_2})
yield
\begin{equation}\label{eq:inner_product_3}
\begin{aligned}
&\frac{\sum_{i=1}^m \langle \lambda_i^kg_i^k, \bar x^k -\theta^*
\rangle}{m}\\
&\geq \frac{\sum_{i=1}^m (-\|x_i^k-\bar{x}^k\|^2
-L^2(\lambda_i^k)^2\|\bar x^k -\theta^*\|^2)}{2m}+\\
&\frac{\sum_{i=1}^m(\langle(\lambda_i^k-\bar\lambda^k)\nabla
f_i(\bar{x}^k), \bar x^k -\theta^*
 \rangle +\langle\bar\lambda^k\nabla f_i(\bar{x}^k), \bar x^k -\theta^*
 \rangle)}{m}\\
 &= -\frac{\sum_{i=1}^m  \|x_i^k-\bar{x}^k\|^2
}{2m}-\frac{L^2\|\lambda^k\|^2\|\bar x^k -\theta^*\|^2 }{2m}+\\
&\frac{\sum_{i=1}^m \langle(\lambda_i^k-\bar\lambda^k)\nabla
f_i(\bar{x}^k), \bar x^k -\theta^*
 \rangle }{m}+\bar\lambda^k\langle\nabla F(\bar x^k), \bar x^k
 -\theta^*\rangle\\
 &\geq -\frac{\sum_{i=1}^m  \|x_i^k-\bar{x}^k\|^2
}{2m}-\frac{L^2\|\lambda^k\|^2\|\bar x^k -\theta^*\|^2 }{2m}+\\
&\frac{\sum_{i=1}^m \langle(\lambda_i^k-\bar\lambda^k)\nabla
f_i(\bar{x}^k), \bar x^k -\theta^*
 \rangle }{m}+\bar\lambda^k(F(\bar x^k)-F(\theta^*)
\end{aligned}
\end{equation}
where we used the convexity of $F(\cdot)$ in the last inequality.

Note  $\lambda^k=[\lambda_1^k,\cdots,\lambda_m^k]^T$, we always have
\begin{equation}\label{eq:inequality_in_Theo1}
\begin{aligned}
&\frac{\sum_{i=1}^m \langle(\lambda_i^k-\bar\lambda^k)\nabla
f_i(\bar{x}^k), \bar x^k -\theta^*
 \rangle }{m}\\
 &=\frac{ \langle\sum_{i=1}^m(\lambda_i^k-\bar\lambda^k)\nabla
f_i(\bar{x}^k), \bar x^k -\theta^*
 \rangle }{m}\\
 &\geq - \frac{ \|\sum_{i=1}^m(\lambda_i^k-\bar\lambda^k)\nabla
f_i(\bar{x}^k)\|\: \|\bar x^k -\theta^*
 \|}{m}\\
 &=- \frac{\| \left((\lambda^k-\bar\lambda^k{\bf 1}_m)\otimes {\bf 1}_d\right)^Tm\nabla
f( {\bf 1}\otimes\bar{x}^k)\|\: \|\bar x^k -\theta^*
 \|}{m}\\
 &\geq - \sqrt{d}\|  \lambda^k-\bar\lambda^k{\bf 1}_m\| \:\|\nabla
f({\bf 1}\otimes\bar{x}^k)\|\: \|\bar x^k -\theta^*
 \|
\end{aligned}
\end{equation}
where $\otimes$ is Kronecker product, and we used $m\nabla
f({\bf 1}\otimes\bar{x}^k))= [(\nabla f_1(\bar{x}^k))^T,\cdots,(\nabla f_m(\bar{x}^k))^T ]^T$ in the last equality and  Cauchy-Schwarz
inequality in the last inequality. Furthermore, $\|\nabla f({\bf
1}\otimes\bar x^k)\|$ can be bounded by using   {   $\nabla
f(x^*)=0$ at $x^*={\bf 1}\otimes\theta^*$} and the Lipschitz
gradient property of $f(\cdot)$, as follows:
\begin{equation}\label{eq:nabla_f}
\begin{aligned}
\|\nabla f({\bf 1}\otimes\bar x^k)\|&= \|\nabla f({\bf 1}\otimes\bar
x^k)-\nabla f(x^*)\|
\\
&\le L\|{\bf 1}\otimes\bar x^k- x^*\|=L\sqrt{m}\, \|\bar x^k -
\theta^*\|
\end{aligned}
\end{equation}

Combining (\ref{eq:inner_product_3}),
(\ref{eq:inequality_in_Theo1}), and (\ref{eq:nabla_f})
 yields
\begin{equation}\label{eq:inner_product_4}
\begin{aligned}
\hspace{-0.2cm}&\frac{\sum_{i=1}^m \langle \lambda_i^kg_i^k, \bar
x^k -\theta^*
\rangle}{m}\\
\hspace{-0.2cm} &\geq -\frac{\sum_{i=1}^m  \|x_i^k-\bar{x}^k\|^2
}{2m}-\frac{L^2\|\lambda^k\|^2\|\bar x^k -\theta^*\|^2 }{2m}-\\
\hspace{-0.2cm}& L\sqrt{md}\|  \lambda^k-\bar\lambda^k{\bf 1}_m\| \:
\|\bar x^k -\theta^*
 \|^2  +\bar\lambda^k(F(\bar x^k)-F(\theta^*)
\end{aligned}
\end{equation}

We next estimate the last term in relation~\eqref{eq-rel1},   for
which we use $\frac{1}{m}g^k=\nabla f(x^k)$, the notation
$\lambda^k=[\lambda_1^k,\cdots,\lambda_m^k]^T$, and the
Cauchy-Schwarz inequality:
\[
\frac{\left\|\sum_{i=1}^m\lambda_i^k g_i^k\right\|^2}{m^2}=  \left\|
(\lambda^k\otimes {\bf 1})^T \nabla f(x^k)\right\|^2 \le
d\|\lambda^k\|^2\|\nabla f(x^k)\|^2
\]
We then add and subtract $\nabla f(x^*)=0$   to obtain
\begin{equation}\label{eq:||g_i||}
\begin{aligned}
\frac{\left\|\sum_{i=1}^m\lambda_i^k g_i^k\right\|^2}{m^2} &\le
2d\|\lambda^k\|^2 \|\nabla f(x^k)-\nabla f(x^*)\|^2 \\
&\le 2d\|\lambda^k\|^2 L^2\|x^k -x^*\|^2
\end{aligned}
\end{equation}
where the last inequality follows by the Lipschitz continuity of
$\nabla f$. Further using the inequality
\begin{equation}\label{eq:x^k-barx_norm}
\begin{aligned}
\|x^k-x^{\ast}\|^2&\leq \|x^k- {\bf 1}\otimes\bar x^k+ {\bf
1}\otimes\bar x^k-x^{\ast}\|^2\\
&\leq 2\|x^k- {\bf 1}\otimes\bar x^k\|^2 +2\|{\bf 1}\otimes\bar
x^k-x^{\ast}\|^2\\
&\leq 2\sum_{i=1}^m \|x_i^k-\bar{x}^k\|^2 +2m\| \bar
x^k-\theta^{\ast}\|^2
\end{aligned}
\end{equation}
we have from (\ref{eq:||g_i||})
\begin{equation}\label{eq:||g_i||_final}
\begin{aligned}
\frac{\left\|\sum_{i=1}^m\lambda_i^k g_i^k\right\|^2}{m^2} \le&
4d\|\lambda^k\|^2 L^2\sum_{i=1}^m \|x_i^k-\bar{x}^k\|^2\\
&+4md\|\lambda^k\|^2 L^2\| \bar x^k-\theta^{\ast}\|^2
\end{aligned}
\end{equation}
Substituting (\ref{eq:inner_product_4}) and (\ref{eq:||g_i||_final})
into  \eqref{eq-rel1}  yields
\begin{equation}\label{eq-rel11}
\begin{aligned}
& \left\|\bar x^{k+1}-\theta^*\right\|^2\leq\left\|\bar x^k
-\theta^* \right\|^2 \\&+\frac{\sum_{i=1}^m \|x_i^k-\bar{x}^k\|^2
}{m}+\frac{L^2\|\lambda^k\|^2\|\bar x^k -\theta^*\|^2 }{m}+\\
& 2L\sqrt{md}\|  \lambda^k-\bar\lambda^k{\bf 1}_m\|  \|\bar x^k
-\theta^*
 \|^2  -2\bar\lambda^k(F(\bar x^k)-F(\theta^*)\\
 &+4d\|\lambda^k\|^2 L^2\sum_{i=1}^m \|x_i^k-\bar{x}^k\|^2 +4md\|\lambda^k\|^2 L^2\| \bar x^k-\theta^{\ast}\|^2
\end{aligned}
\end{equation}

We can group the common terms on the right hand side of the
preceding relation and obtain
\begin{equation}\label{eq:DS_Step1_final}
\begin{aligned}
& \left\|\bar x^{k+1}-\theta^*\right\|^2\leq\left\|\bar x^k -\theta^* \right\|^2\times\\
& \left(1+ \frac{L^2\|\lambda^k\|^2+2Lm\sqrt{md}\|
\lambda^k-\bar\lambda^k{\bf 1}\|}{m}+4md\|\lambda^k\|^2 L^2\right)
\\&+(\frac{1}{m}+4d\|\lambda^k\|^2 L^2) \sum_{i=1}^m\|x_i^k-\bar{x}^k\|^2
 -2\bar\lambda^k(F(\bar x^k)-F(\theta^*)
\end{aligned}
\end{equation}

Step II: Relationship for $\sum_{i=1}^m\|x_i^k-\bar x^k\|^2$. For
the convenience of analysis, we write
PDG-DS on per-coordinate expressions. Define for all
$\ell=1,\ldots,d,$ and $k\ge0$,
\[
\begin{aligned}
x^k(\ell)&=\left[[x_1^k]_\ell,\ldots,[x_m^k]_\ell\right]^T,\:
g^k(\ell)&=\left[[g_1^k]_\ell,\ldots,[g_m^k]_\ell\right]^T
\end{aligned}
\] In
this per-coordinate view, (\ref{eq:decaying_stepsize}) and
(\ref{eq:bar_x^k}) have the following form  for all
$\ell=1,\ldots,d,$ and $k\ge0$:
\begin{equation}\label{percord_DS}
\begin{aligned}
x^{k+1}(\ell)&=Wx^k(\ell)-B^k \Lambda^{k} g^{k}(\ell)\\
[\bar x^{k+1}]_{\ell}&=[\bar x^k]_{\ell} -\frac{1}{m}{\bf
1}^T\Lambda^kg^k(\ell)
\end{aligned}
\end{equation}
From (\ref{percord_DS}), we obtain
\[
\begin{aligned}
x^{k+1}(\ell)-[\bar x^{k+1}]_{\ell}{\bf 1}=&Wx^k(\ell)-[\bar
x^k]_{\ell} {\bf 1} \\
&-\left(B^k\Lambda^k g^k(\ell) -\frac{1}{m}{\bf
1}^T\Lambda^kg^k(\ell){\bf1} \right)
\end{aligned}
\]

Noting that $[\bar x^k]_{\ell} {\bf 1}=\frac{1}{m}{\bf 1}{\bf
1}^Tx^k(\ell)$ and $ \frac{1}{m}{\bf 1}^T\Lambda^k g^k(\ell) {\bf1}=
\frac{1}{m}{\bf1}{\bf 1}^T\Lambda^k g^k(\ell)$, we have
\[
\begin{aligned}
x^{k+1}(\ell)-[\bar x^{k+1}]_{\ell}{\bf 1}=&\bar Wx^k(\ell) -\bar
B^k\Lambda^k g^k(\ell)
\end{aligned}
\]
where $\bar W=W- \frac{ {\bf 1}{\bf 1}^T}{m}$ and $\bar B^k=B^k-
\frac{{\bf1}{\bf 1}^T }{m}$.

Noticing   $\bar W[\bar x^k]_{\ell}{\bf1}=\left(W- \frac{1}{m}{\bf
1}{\bf 1}^T\right) [\bar x^k]_{\ell}{\bf1}=0$, by subtracting this
expression from the right hand side of the preceding relation, we
obtain
\[
\begin{aligned}
x^{k+1}(\ell)-[\bar x^{k+1}]_{\ell}{\bf 1}=&\bar
W(x^k(\ell)-[\bar{x}^k]_{\ell}{\bf 1}) -\bar B^k\Lambda^k g^k(\ell)
\end{aligned}
\]

Taking  norm on both sides and using   $\eta=\|W- \frac{1}{m}{\bf
1}{\bf 1}^T\|$ from Assumption \ref{assumption:W}, we obtain
\begin{equation}\label{eq:B^k_algorithm 1}
\begin{aligned}
\|x^{k+1}\hspace{-0.07cm}(\ell)\hspace{-0.1cm}-\hspace{-0.1cm}[\bar x^{k+1}]_{\ell}{\bf 1}\|\le&
\eta\|x^k(\ell)\hspace{-0.1cm}-\hspace{-0.1cm}[\bar x^k]_{\ell}{\bf1}\|
\hspace{-0.1cm}+\hspace{-0.1cm}\left\|\bar B^k\right\| \|\Lambda^k\|\, \|g^k(\ell)\|
\end{aligned}
\end{equation}
The column stochastic property of $B^k$ implies
\begin{equation}\label{eq:B^k_norm}
\left\|\bar B^k\right\|\le \left\|\bar B^k\right\|_F
 \leq m
\end{equation}
where $\|\cdot\|_F$ denotes the Frobenius matrix norm, yielding
\[
\begin{aligned}
\|x^{k+1}(\ell)-[\bar x^{k+1}]_{\ell}{\bf 1}\|\le&
\eta\|x^k(\ell)-[\bar x^k]_{\ell}{\bf1}\|
 +m\|\Lambda^k\|\, \|g^k(\ell)\|
\end{aligned}
\]

By taking squares on both sides  and using the inequality $2ab\le
\epsilon a^2 + \epsilon^{-1}b^2$ valid for any  $a$, $b$, and
$\epsilon>0$, we obtain
\[
\begin{aligned}
\|x^{k+1}(\ell)-[\bar x^{k+1}]_{\ell}{\bf 1}\|^2\le&
\eta^2(1+\epsilon)\|x^k(\ell)-[\bar x^k]_{\ell}{\bf1}\|^2\\
& +m^2(1+\epsilon^{-1})\|\Lambda^k\|^2\|g^k(\ell)\|^2
\end{aligned}
\]

Summing these relations over $\ell=1,\ldots,d$, and noting  $
\sum_{\ell=1}^d \|x^k(\ell) - [\bar x^k]_\ell{\bf
1}\|^2=\sum_{i=1}^m\|x^k_i - \bar x^k \|^2$ and $\sum_{\ell=1}^d
\|g^k(\ell)\|^2=\sum_{i=1}^m\|g^k_i\|^2$, we obtain
\begin{equation}\label{eq:x^{k+1}-barx_norm}
\begin{aligned}
\sum_{i=1}^m\|x^{k+1}_i - \bar x^{k+1} \|^2\le&
\eta^2(1+\epsilon)\sum_{i=1}^m\|x^k_i - \bar x^k \|^2 \\
& +m^2(1+\epsilon^{-1})\|\Lambda^k\|^2\|g^k\|^2
\end{aligned}
\end{equation}

We next focus on estimating $\|g^k\|^2$. Noting that $g^k=m\nabla
f(x^k)$, $\nabla f(x^*)=0$, and $f$ has Lipschitz continuous
gradients, we have
\begin{equation}\label{eq:g^k_norm}
\begin{aligned}
\|g^k\|^2 &=m^2\|\nabla f(x^k)-\nabla f(x^*)\|^2 \le m^2L^2\|x^k-
x^*\|^2\\
&\leq 2m^2L^2 \sum_{i=1}^m \|x_i^k-\bar{x}^k\|^2 +2m^3L^2\| \bar
x^k-\theta^{\ast}\|^2
\end{aligned}
\end{equation}
where the last inequality used (\ref{eq:x^k-barx_norm}).

Substituting (\ref{eq:g^k_norm}) into (\ref{eq:x^{k+1}-barx_norm})
and grouping terms yield
\[
\begin{aligned}
&\sum_{i=1}^m\|x^{k+1}_i - \bar x^{k+1}
\|^2\le2m^5L^2(1+\epsilon^{-1})\|\Lambda^k\|^2\| \bar
x^k-\theta^{\ast}\|^2\\
&+\left(\eta^2(1+\epsilon)+2m^4L^2(1+\epsilon^{-1})\|\Lambda^k\|^2)\right)\sum_{i=1}^m\|x^k_i
- \bar x^k \|^2
\end{aligned}
\]

By letting $\epsilon=\frac{1-\eta}{\eta}$ with $\epsilon>0$, and
noting $\eta\in(0,1)$, $1+\epsilon=\eta^{-1}$,
$1+\epsilon^{-1}=(1-\eta)^{-1}$, and
$\|\Lambda^k\|=\max_i\lambda_i^k\le \|\lambda^k\|$, we arrive at
\begin{equation}\label{eq:DS_Step2_final}
\begin{aligned}
&\sum_{i=1}^m\|x^{k+1}_i - \bar x^{k+1}
\|^2\le2m^5L^2(1-\eta)^{-1}\|\lambda^k\|^2\| \bar
x^k-\theta^{\ast}\|^2\\
&\qquad+\left(\eta
+2m^4L^2(1-\eta)^{-1}\|\lambda^k\|^2)\right)\sum_{i=1}^m\|x^k_i -
\bar x^k \|^2
\end{aligned}
\end{equation}

Combining (\ref{eq:DS_Step1_final}) and (\ref{eq:DS_Step2_final}),
we have
\begin{equation}\label{algorithm_proof_final}
\bv^{k+1} \le \left( \left[\begin{array}{cc} 1 & \frac{1}{m}\cr 0&
\eta\cr
\end{array}\right]
+A^k \right)\bv^k  - 2\bar\lambda^k\left[\begin{array}{c}
 (F(\bar x^k)-F(\theta^*)\cr 0\end{array}\right]
\end{equation}
where $\bv^k=\left[\begin{array}{c}\|\bar x^k-\theta^*\|^2\cr
\sum_{i=1}^m\|x_i^k-\bar x^k\|^2\end{array}\right]$, $ A^k=
\left[\begin{array}{cc}A_{11}&A_{12}\\A_{21}&A_{22}\end{array}\right]
$ with $A_{11}\triangleq\frac{L^2\|\lambda^k\|^2+2Lm\sqrt{md}\|
\lambda^k-\bar\lambda^k{\bf 1}\|}{m}+4md\|\lambda^k\|^2 L^2$,
$A_{12}=4d\|\lambda^k\|^2 L^2$,
$A_{21}=2m^5L^2(1-\eta)^{-1}\|\lambda^k\|^2$, and
$A_{22}=2m^4L^2(1-\eta)^{-1}\|\lambda^k\|^2$.

Because $\left[\begin{array}{cc} 1 & \frac{1}{m}\cr 0& \eta\cr
\end{array}\right]\leq \left[\begin{array}{cc} 1 & 1\cr 0&
\eta\cr
\end{array}\right]$ and $A^k\leq a^k {\bf 1}{\bf 1}^T $ hold when
$a^k$ is set to
$a^k=\max\left\{A_{11},A_{12},A_{21},A_{22}\right\}$, we can see
that (\ref{eq-fin}) in Proposition \ref{th-main} is satisfied. Further
note that under the conditions in the statement, all conditions for
$\{a^k\}$, $\{b^k\}$, and $\{c^k\}$ in Proposition \ref{th-main} are
also satisfied ($b^k$ is always $0$ here). Therefore, we have the
claimed results.
\end{proof}

{
}

{
}
\begin{Remark 1}
To our knowledge, for decentralized gradient methods with
diminishing stepsizes, our result is the first to prove exact
convergence  under general time-varying stepsize heterogeneity. In
fact, the condition in (\ref{eq:stepsize_heterogenerity}) can be
satisfied even when the stepsize differences in a finite number of
iterations are arbitrarily large. This can enable strong privacy, as
detailed in Sec. \ref{se:privacy_DS}.
\end{Remark 1}

\subsection{Privacy Analysis}\label{se:privacy_DS}
Recall that in Sec. II we identify the gradients of  agents as
information to be protected in decentralized optimization. In this
subsection, we will show that the PDG-DS algorithm can effectively
protect the gradients of all participating agents from being
inferable by  honest-but-curious adversaries and external
eavesdroppers. To this end, we first give a privacy metric and our
definition of privacy protection.

{  Our result in Theorem \ref{theorem_PDGD_DS_convergence} ensures deterministic convergence of all agents to an exact optimal solution, without imposing any stochastic conditions (e.g., types of stochastic distribution or population parameters such as mean and variance) on the uncertainty in stepsizes. To be consistent with such a deterministic analysis framework, inspired by the logarithmic difference
metric used in differential privacy, we propose a deterministic framework for privacy analysis. } More specifically, we define the difference
{ from $x$ to $x'$ as follows (in the log scale, could be positive or negative}):
\begin{equation}\label{eq:difference}
\zeta= \log\textstyle\frac{\|x\|}{\|x'\|}
\end{equation}

\begin{Definition 1}\label{definition:privacy}
For a network of $m$ agents in decentralized optimization, the
privacy of agent $i$ is preserved if for { any finite number of iterations $T$,
   its gradient values $g_i^1,\cdots, g_i^T$    always have
  alternative realizations  $\hat{g}_i^1,\cdots, \hat{g}_i^T$}   which allow each $\hat{g}_i^k$ ($1\leq k\leq T$) to  have an arbitrarily large difference {(under the log scale defined in (\ref{eq:difference}), could be different for different $k$)} from $g_i^k$, but lead to the
same shared information in inter-agent communications.
\end{Definition 1}

The above privacy definition requires that when an agent's gradient
is perturbed by an arbitrary value $\zeta$, its shared information
can still be the same, i.e., an alteration to an agent's gradient
value is not distinguishable by an adversary having access to all
information shared by the agent. Since the alteration in gradient
can be arbitrarily large, our privacy definition   requires that an
adversary cannot even find a range for a protected value, and hence
is more stringent than many existing privacy definitions (e.g.,
\cite{zhang2019admm,gade2018private}) that only require an adversary
unable to {\it uniquely} determine a protected value.

\begin{Theorem 1}\label{theo:privacy_diminishing}
In the presence of honest-but-curious or eavesdropping adversaries,
  PDG-DS can protect the
privacy of all participating agents defined in Definition
\ref{definition:privacy}.
\end{Theorem 1}

\begin{proof}
Without loss of generality, we first consider the protection of the gradient of agent $i$ at any {single} time instant $k$, and then { show that the argument also applies to any finite number of time instants (iterations)}.  When the gradient  is $g^k_i$, we
represent the information that agent $i$ shares with neighboring
agents when participating  PDG-DS as $\mathcal{I}_i$. According to Definition
\ref{definition:privacy}, we  have to prove that when the gradient
is altered { from $g_i^k$} to $\hat{g}^k_i=e^{ \zeta^k}g^k_i$ with $\zeta^k$ difference from $g^k_i$ according to the log-scale metric in (\ref{eq:difference}),
the corresponding shared information $\hat{\mathcal{I}}_i$ of agent
$i$ when running the algorithm
 could be identical to $\mathcal{I}_i$ under
any $\zeta^k>0$.

According to Algorithm PDG-DS, agent $i$ shares the following
information in decentralized optimization:
\[
\mathcal{I}_i=\mathcal{I}_i^{\rm sent}\bigcup\mathcal{I}_i^{\rm
public}
\]
with $ \mathcal{I}_i^{\rm sent} =\left\{v^k_{ji}\triangleq
w^k_{ji}x_i^k-b_{ji}^k\lambda_i^kg_i^k|k=1,2,\cdots\right\} $ and $
\mathcal{I}_i^{\rm public}=\left\{W\bigcup
\sum_{j\in\mathbb{N}_i}b^k_{ji}=1|k=0,1,\cdots\right\} $.  One can
obtain that  at some iteration $k$, if the gradient is changed to
$\hat{g}^k_i=e^{ \zeta^k}g^k_i$, the difference defined in
(\ref{eq:difference}) is $ \zeta^k$. However, in this case, if we set
the stepsize $\hat{\lambda}_i^k$ to
$\hat{\lambda}_i^k=e^{ -\zeta^k}\lambda_i^k$, then the corresponding
shared information will still be $v_{ji}^k$. Since other parameters
are not changed and changing the stepsize {  from $\lambda_i^k$} to
$\hat{\lambda}_i^k=e^{ -\zeta^k}\lambda_i^k$ will not violate the
summable stepsize heterogeneity condition in
(\ref{eq:stepsize_heterogenerity}) for any given $ \zeta^k<\infty$,
according to Theorem \ref{theorem_PDGD_DS_convergence}, convergence
to the optimal solution will still be guaranteed. Therefore, changes
in an agent's gradient   can be completely covered by the agent's
flexibility in changing its stepsize, which does not affect the
convergence. Thus, privacy of any agent's   gradient will be protect
when running PDG-DS. { Given that the
summable stepsize heterogeneity condition in
(\ref{eq:stepsize_heterogenerity}) allows the stepsize of agent $i$ to change by any finite amount for any finite number of iterations
according to Theorem \ref{theorem_PDGD_DS_convergence}, one can obtain that the  privacy of every agent's gradients in any number of iterations can be completely covered by the flexibility in changing the agent's stepsize in these iterations, as long as the  number of these iterations is finite.}
{It is worth noting that the perturbation does not violate the convexity and Lipschitz conditions in Assumption 2. This is because in order for an adversary to check if Assumption 2 is violated, it has to know $x_i^k$ and $L$, which, however, are not available to adversaries: before convergence, $x_i^k$ is inaccessible to the adversary because the information shared by agent $i$ is $w_{ji}x_i^k-b_{ji}^k\lambda_i^kg_i^k$, avoiding $x_i^k$ from being inferable; $L$ is inaccessible  to the adversary either because Assumption 2 only requires  all $f_i$ to have  finite Lipschitz constants, and agents do not share their Lipschitz constants (agents do not need to know the upper bound on Lipschitz constant $L$ in Assumption 2) in the implementation of the algorithm.  In fact, even with the gradient $g_i^k$ unchanged, the value of observation    $w^k_{ji}x_i^k-b_{ji}^k\lambda_i^kg_i^k$ in Algorithm 1   can be changed by an arbitrary finite value by changing the stepsize $\lambda_i^k$.
 Therefore, before convergence, an adversary cannot use Assumption 2 to confine the change in observed values and further confine the change in the value of $g_i^k$.  After convergence, the perturbation does not violate the convexity and Lipschitz conditions, either.  In fact, although $x_i^k$ becomes accessible to the adversary after convergence,   gradient information  is eliminated in adversary's observation (the shared information $w_{ji}x_i^k-b_{ji}^k\lambda_i^kg_i^k$   becomes $w_{ji}x_i^k$) because $\lambda_i^k$ converges to zero. So after convergence,  the adversary still cannot use Assumption 2 to confine  changes in gradients.}
\end{proof}

\begin{Remark 1}\label{re:convergence_privacy}
Even after convergence when $g_i^k$ becomes a constant, an adversary still cannot infer gradients from shared messages in PDG-DS. More specifically,  when $g_i^k$ converges to a constant value, the stepsize $\lambda_i^k$  also converges to zero, which completely eliminates the information of $g_i^k$  in observed information (the observed information becomes $w_{ji}x_i^k$ after convergence). This can also be understood intuitively as follows: Even if the adversary can collect  $T\rightarrow\infty$  observations $w_{ji}x_i^k-b_{ji}^k\lambda_i^kg_i^k$ in the neighborhood of the optimal point and establish  a system  of $T$ equations to solve for $g_i^k$ (which can be  viewed  to be   approximately time-invariant  in the neighborhood of the optimal point), the number of unknowns $b_{ji}^k$, $\lambda_i^k$, and $g_i^k$ in the system  of $T$ equations is $3T$ (even if we view $\lambda_i^k$ and $g_i^k$ approximately as   constants in the neighborhood of the optimal point, the number of unknowns is still $T+2$), which makes it impossible for the adversary to solve for $g_i^k$ using the  system  of $T$  equations established from observations.
\end{Remark 1}

{
\begin{Remark 1}
  Different from existing privacy solutions for decentralized optimization that patch  a privacy mechanism (e.g., differential-privacy noise or encryption) with a  pre-designed decentralized  optimization algorithm, our proposed algorithm uses stepsize and coupling coefficients that are inherent to the decentralized optimization algorithm to  perturb gradients, and hence has inherent privacy.
\end{Remark 1}

}
{
\begin{Remark 1}
  Different from  differential privacy which uses some index to    quantitatively measure the level of privacy, our defined privacy is   binary in the sense that it can only be satisfied or not. It is also deterministic,  employing uncertainties in random  stepsize and $B^k$ to establish obscuration.  Such binary and deterministic definition of privacy   has also been widely used in other privacy frameworks such as $\ell$-diversity \cite{machanavajjhala2007diversity}  and Shamir's secret sharing \cite{shamir1979share}, which also employ random variables to establish obscuration.  Note that our deterministic privacy framework allows us to establish deterministic convergence to the exact optimal solution, which  is in stark contrast to differential privacy that inevitably compromises the accuracy of optimization. In fact, for iterative distributed optimization and learning algorithms, differential privacy has to consider  the cumulative privacy loss in all iterations, which usually leads to overly conservative design and makes it  sometimes regarded as too restrictive  \cite{melis2019exploiting,jayaraman2019evaluating}.
\end{Remark 1}
}

\begin{Remark 1}
Existing accuracy-maintaining privacy  approaches for decentralized
optimization can only protect the privacy of participating agents
when  the interaction topology meets certain conditions. For
example, the  approach in \cite{gade2018private} assumes that an
adversary cannot have access to messages sent on at least one
communication channel of an agent to guarantee the privacy of this
agent. The  approach in \cite{zhang2018enabling} requires that an
adversary cannot be the only neighbor of a target agent. To the
contrary, our PDG-DS  can protect the privacy of an agent without
any constraint on the interaction topology. In fact, to our
knowledge, our algorithm is the first decentralized gradient based
algorithm that can guarantee both optimization accuracy and  privacy
defined in Definition \ref{definition:privacy}
 when an adversary  has access to all
shared information.
\end{Remark 1}
\begin{Remark 1}
In line with the discussions in the previous remark, PDG-DS can
protect the defined privacy of participating agents irrespective of
the number of adversaries and existence of collusion among
adversaries.
\end{Remark 1}

{
\begin{Remark 1}
Our privacy approach protects the gradient value by making an arbitrarily large perturbation on it indistinguishable to  adversaries having   access to all information shared by the agent. This is important  in applications where the gradient value carries sensitive information. For example, in the rendezvous problem, the objective function takes the form $f_i(x)=\|x_{i,0}-x\|^2$ with $x_{i,0}\in\mathbb{R}^2$ denoting the initial position (which is sensitive and has to be protected) and $x$ denoting the optimization variable \cite{ruan2019secure,Huang15}.  It can be seen that an adversary can easily calculate the private position $x_{i,0}$ if it has access to the gradient value $g_i=2(x-x_{i,0})$ and optimization variable $x$ (note that $x$ is  disclosed explicitly in most existing distributed optimization algorithms). The same argument holds for sensor-network based localization applications where disclosing the gradient value leads to disclosing the positions of sensors (which are sensitive and should be protected) \cite{zhang2018enabling}. In scenarios where gradient directions instead of gradient values carry sensitive information, we may have to resort to other privacy approaches such as \cite{Huang15}.
\end{Remark 1}
}

\section{An inherently privacy-preserving decentralized gradient  algorithm with non-diminishing stepsizes}

Because  diminishing stepsizes in decentralized gradient methods may
slow down convergence, plenty of efforts have been devoted to
developing decentralized optimization algorithms that can achieve
accurate optimization results under a non-diminishing stepsize.
Typical examples include  gradient-tracking based algorithms such as
Aug-DGM \cite{xu2015augmented}, DIGing
\cite{nedic2017achieving}, AsynDGM
\cite{xu2017convergence}, AB \cite{xin2018linear}, Push-Pull
\cite{pu2020push,du2018accelerated}, and ADD-OPT \cite{xi2017add},
etc. However, these algorithms will lead to privacy breaches in
implementation. For example, in DIGing \cite{nedic2017achieving},
agent $i$ employs two variables (the optimization variable $x_i^k$
and the auxiliary gradient-tracking variable $y^k_i$) implementing
the following update rule (note that under our assumption, $i \in
\mathbb{N}_i$):
 \[
\left\{
\begin{aligned}
x_i^{k+1}&= \sum\nolimits_{j\in\mathbb{N}_i}w_{ij}x_j^k-\lambda
y_i^k\\y_i^{k+1}&=
\sum\nolimits_{j\in\mathbb{N}_i}w_{ij}(y_j^k+g_j^{k+1}-g_j^k)
\end{aligned}
\right.
\]
At iteration $k=0$, agent $j$ sets $y_j^0=g_j^0$ and sends
$w_{ij}(y_i^0+g_j^1-g_j^0)=w_{ij}g_j^1$ to its neighboring agent
$i$. At iteration $k=1$, agent $j$ further sends  $x_j^1$ to agent
$i$. Given that $w_{ij}$s are publicly known, agent $i$ can easily
determine the gradient of agent $j$ at $x_j^1$. Using a similar
argument, we can  see that other commonly used gradient-tracking
based algorithms also have the same issue of leaking  agents'
gradient information, even when the stepsizes are heterogeneous (see
Sec. IV.B for details).

 The EXTRA algorithm \cite{shi2015extra} can also ensure
 convergence to the exact optimal solution under non-diminishing
 stepsizes:
 \[
 x_i^{k+2}=\sum\nolimits_{j\in\mathbb{N}_i}w_{1,ij}x_j^{k+1}-\sum\nolimits_{j\in\mathbb{N}_i}w_{2,ij}x_j^k-\lambda(g_i^{k+1}-g_i^{k})
 \]
However, since  $\lambda,w_{1,ij},w_{2,ij}$ are publicly known, and
an agent $i$ has to share $x_i^k$ directly, one can see that the
gradient information of participating agent $i$ will also be
disclosed to neighboring agents and eavesdroppers.

 Motivated by the observation that the main sources of information
 leakage  in decentralized optimization are constant parameters and
 the sharing of
two messages by every agent  in every iteration, we propose the
following inherently privacy-preserving decentralized gradient based
algorithm which can protect the gradients of participating agents
while ensuring convergence to the exact optimal solution under
non-diminishing stepsizes (the per-agent version is given in
Algorithm PDG-NDS):
\begin{equation}\label{eq:nondiminishing_stepsize}
\begin{aligned}
&x^{k+2}=2(W\otimes I_d)x^{k+1}-(W^2\otimes
I_d)x^k\\
&\qquad-\left(\left((B^k\Lambda^{k+1})\otimes
I_d\right)g^{k+1}-\left((B^k\Lambda^{k})\otimes
I_d\right)g^{k}\right)
\end{aligned}
\end{equation}
where $B^k=\{b^k_{ij}\}\in\mathbb{R}^{m\times m}$ is a
column-stochastic matrix, $\Lambda^k={\rm
diag}[\lambda_1^k,\lambda_2^k,\cdots,\lambda_m^k]$ with
$\lambda_i^k\geq 0$ denoting the stepsize of agent $i$ at iteration
$k$, $g^k=[(g_1^k)^T,(g_2^k)^T,\cdots (g_m^k)^T]^T$, and
$\otimes$ denotes Kronecker product.

\noindent\rule{0.49\textwidth}{0.5pt} \noindent\textbf{PDG-NDS:
Privacy-preserving decentralized gradient method with
non-diminishing stepsizes}

\vspace{-0.2cm}\noindent\rule{0.49\textwidth}{0.5pt}
\begin{enumerate}
    \item[] Public parameters: $W$
    \item[] Private parameters for each agent $i$:  $b^k_{ji}$, $\lambda_i^k$, and {$x_i^0$}
    \item At iteration $k=1$: Agent $i$ shares $x_i^0$ (randomly selected) with neighbors and updates its state as follows

     \[
     x^1_i=\sum\nolimits_{j\in\mathbb{N}^i}w_{ij}x_j^0-\lambda_i^0\nabla f_i(x_i^0)
     \]
    \item {\bf for  $k=2,3,\cdots$ do}
    \begin{enumerate}
        \item Every agent $j$ computes and sends $v^k_{ij}$  (defined in (\ref{eq:v_ij}))
        to
        all
        agents
        $i\in\mathbb{N}_j$ where  $\{W^2\}_{ij}$ denotes the $(i,j)$th entry of matrix $W^2$:
        \begin{equation}\label{eq:v_ij}
        \begin{aligned}
        v^k_{ij}&=2w_{ij}x_j^{k-1}+\{W^2\}_{ij}x_j^{k-2}\\
        &\qquad -b_{ij}^{k-2}(\lambda_j^{k-1}g_j^{k-1}-\lambda_j^{k-2}g_j^{k-2})
        \end{aligned}
        \end{equation}
        \item After receiving $v^k_{ij}$ from all $j\in\mathbb{N}_i$, agent $i$ updates its state as follows:
        \begin{equation}\label{update-rule-01}
            x_i^k=\sum\nolimits_{j\in\mathbb{N}_i}v^k_{ij}
        \end{equation}
        \vspace{-0.5cm}
        \item {\bf end}
    \end{enumerate}
\end{enumerate}
\vspace{-0.2cm}\rule{0.49\textwidth}{0.5pt}

\begin{Remark 1}
When $B^k$ is set to $I_m$ and  stepsizes are the same and constant,
i.e., $\lambda_i^k=\lambda$, PDG-NDS reduces to EXTRA
\cite{shi2015extra}.
\end{Remark 1}
\begin{Remark 1}
PDG-NDS  requires one agent to share only one variable with every
neighboring agent at each iteration. This  is different from all
existing gradient-tracking based algorithms  which have to exchange
two variables  between two neighboring agents in every iteration
(one optimization variable and one auxiliary variable tracking the
gradient of the global objective function). This difference is key
to 1) reduce communication overhead; 2) enable privacy because
exchanging the additional gradient-tracking variable will disclose
gradient information, as detailed in Sec. IV.B.
\end{Remark 1}
\subsection{Convergence analysis}
We define an auxiliary variable
\begin{equation}\label{eq:y_k_definition}
y^k\triangleq (W\otimes I_d) x^{k}-x^{k+1}
\end{equation}

It can be verified that
\begin{equation}\label{eq:y_k+1}
\begin{aligned}
y^{k+1}&=(W\otimes I_d)y^k+\\
&\left(\left((B^k\Lambda^{k+1})\otimes
I_d\right)g^{k+1}-\left((B^k\Lambda^{k})\otimes
I_d\right)g^{k}\right)
\end{aligned}
\end{equation}

Define  mean vectors of $x_i^k$ and $y_i^k$ as $\bar
x^k=\frac{1}{m}\sum_{i=1}^m x_i^k$ and $\bar
y^k=\frac{1}{m}\sum_{i=1}^m y_i^k$, respectively. Then from
(\ref{eq:y_k_definition}), we have
\begin{equation}\label{eq:bar_y_k}
\bar{x}^{k+1}=\bar{x}^k-\bar{y}^k
\end{equation}
 Further using  (\ref{eq:y_k+1}) and the initialization condition
 $x^1=Wx^0-\Lambda^0g^0$ in PDG-NDS,  we have $\bar{y}^0=
 \frac{1}{m}\sum_{i=1}^m\lambda_i^0g_i^0$ and
\begin{equation}\label{eq:eq-bary}
\bar{y}^k=
\textstyle\frac{1}{m}\sum\nolimits_{i=1}^m\lambda_i^kg_i^k
\end{equation}

To prove convergence of our algorithm, we first present two lemmas
and one theorem. The theorem  applies to general distributed
algorithms for solving optimization problem
(\ref{eq:optimization_formulation1}).

\begin{Lemma 1}\label{th-dsystem}
Let $\{\bv^k\}\subset \mathbb{R}^d$ and $\{\bu^k\}\subset
\mathbb{R}^p$ be sequences of non-negative vectors such that
\begin{equation}\label{eq:v^k+1}
\bv^{k+1}\le (V^k+a^k{\bf 1}{\bf1}^T)\bv^k +b^k{\bf 1} -C^k\bu^k,
\forall k\ge0
\end{equation}
where $\{V^k\}$ is a sequence of non-negative matrices, and
$\{a^k\}$ and $\{b^k\}$ are non-negative scalar sequences satisfying
$\sum_{k=0}^\infty  a^k<\infty$ and $\sum_{k=0}^\infty b^k<\infty$.
Assume that there exists a vector $\pi>0$ such that $\pi^T V^k\le
\pi^T$ and $\pi^TC^k \ge 0$ hold for all $k\ge0$. Then,
$\lim_{k\to\infty}\pi^T\bv^k$ exists, the sequence $\{\bv^k\}$ is
bounded, and $\sum_{k=1}^\infty \pi^TC^k\bu^k<\infty$.
\end{Lemma 1}

\begin{proof}
By multiplying (\ref{eq:v^k+1}) with $\pi^T$ and using the
assumptions $\pi^TV^k\le\pi^T$ and $\bv^k\geq 0$, we obtain for
$\forall k\ge0$
\[\pi^T\bv^{k+1}\le \pi^T\bv^k +a^k(\pi^T {\bf 1})({\bf1}^T\bv^k) +b^k\pi^T{\bf 1}-\pi^TC^k\bu^k\] Since $\pi>0$, we have $\pi_{\min}=\min_i\pi_i>0$, and it
follows

\[
{\bf1}^T\bv^k=\frac{1}{\pi_{\min}}
\,\pi_{\min}{\bf1}^T\bv^k\le \frac{1}{\pi_{\min}} \,\pi^T\bv^k
\]
where the inequality in the preceding relation holds since
$\bv^k\ge0$. Therefore, \be\label{eq-oo} \pi^T\bv^{k+1}\le
\left(1+a^k\frac{\pi^T {\bf 1}}{\pi_{\min}}\right)\pi^T\bv^k
+b^k\pi^T{\bf 1} -\pi^TC^k\bu^k, \:\forall k\ge0\ee
 By our assumption,  $\pi^TC^k\bu^k\ge0$ for
all $k$, so (\ref{eq-oo}) implies that the conditions of
Lemma~\ref{lem-polyak} in the Appendix are satisfied with $v^k =
\pi^T\bv^k$, $\a^k = a^k\pi^T {\bf 1}/\pi_{\min}$, and $\b^k = b^k
\pi^T{\bf 1}$. Thus, by Lemma~\ref{lem-polyak}, it follows that
 $\lim_{k\to\infty}\pi^Tv^k$ exists. Consequently, $\{\pi^Tv^k\}$ is
bounded, and under $\pi>0$, implying that $\{\bv^k\}$ is also
bounded. Moreover, by summing the relations in~\eqref{eq-oo}, we
find   $\sum_{k=1}^\infty \pi^TC^k\bu^k<\infty$.
\end{proof}

\begin{Lemma 1}\label{th-dsystem2}
Let $\{\bv^k\}\subset \mathbb{R}^d$  be a sequence of non-negative
vectors such that for $\forall k\ge 0$
\begin{equation}\label{eq:v^k+1_Lemma2}
\bv^{k+1}\le V^k\bv^k +b^k{\bf 1}
\end{equation}
where $\{V^k\}$ is a sequence of non-negative matrices.
Assume that there exist  a vector $\pi>0$ and a scalar sequence
$\{\a^k\}$ such that  $\a^k\in(0,1)$, $\sum_{k=0}^\infty
\a^k=+\infty$, $\lim_{k\to\infty}b^k/\a^k=0$, and $\pi^T V^k\le
(1-\a^k) \pi^T$ for all $k\ge0$. Then, $\lim_{k\to\infty}\bv^k=0$.
\end{Lemma 1}
\begin{proof}
We use Lemma~\ref{lem-polyak2} in the Appendix to establish the
result. By multiplying (\ref{eq:v^k+1_Lemma2}) with $\pi^T$ and
using the assumptions $\pi^TV^k\le(1-\a^k)\pi^T$ and $\bv^k\geq 0$,
we obtain
\[\pi^T\bv^{k+1}\le (1-\a^k)\pi^T\bv^k +b^k\pi^T{\bf 1}, \forall k\ge0\]
Since $\a^k\in(0,1)$, $\sum_{k=0}^\infty \a^k=+\infty$,
$\lim_{k\to\infty}b^k/\a^k=0$, the conditions of
Lemma~\ref{lem-polyak2} in the Appendix are satisfied with $v^k =
\pi^T\bv^k$ and $\b^k = b^k \pi^T{\bf 1}$. Thus,  it follows
$\lim_{k\to\infty}\pi^T\bv^k=0$ and further (because $\pi>0$)
$\lim_{k\to\infty} \bv^k=0$ .
\end{proof}

\begin{Proposition 1}\label{th-main2}
Assume that problem (1) has an optimal solution and that  $F(\cdot)$
in (\ref{eq:optimization_formulation1}) is continuously
differentiable. Suppose that a distributed algorithm generates
sequences $\{x_i^k\}\subseteq\mathbb{R}^d$ and
$\{y_i^k\}\subseteq\mathbb{R}^d$ such that the following relation is
satisfied for any optimal solution $\theta^*$ and for all $k\ge0$,
\begin{equation}\label{eq-fine}
\bv^{k+1} \le \left(V + a^k {\bf 1}{\bf 1}^T\right)\bv^{k}+b^k{\bf
1} - C \left[\begin{array}{c}
 \|\nabla F(\bar x^k)\|^2\cr
 \|\bar y^k\|^2\end{array}\right]
\end{equation}
where $\nu>0$,
\[
\begin{aligned}
\bv^k&\triangleq\left[\hspace{-0.2cm}\begin{array}{c} \nu(F(\bar
x^{k})-F(\theta^*))\cr \sum_{i=1}^m\|x_i^{k}-\bar x^{k}\|^2\cr
\sum_{i=1}^m\|y_i^{k}-\bar y^{k}\|^2\end{array}
\hspace{-0.2cm}\right],
 C=\left[\begin{array}{cc} \gamma^k & \frac{1-\tau^k}{\tau^k}\cr 0 &
0\cr 0 & -(1-\eta)^2(1-c)
\end{array}\hspace{-0.2cm}\right]
\\
&V=\left[\begin{array}{ccc} 1 & 1-\eta &0\cr 0 & \eta
&\frac{1}{1-\eta}\cr 0 & (1-\eta)^2(1-c)(1-\delta) &
c\end{array}\right]
\end{aligned}
\] with $\eta,c,\delta\in(0,1)$, while the scalar
sequences $\{a^k\}$, $\{b^k\}$, $\{\tau^k\}$, $\{\gamma^k\}$ are
nonnegative  satisfying $\tau^k\in(0,1)$,
 $\frac{1-\tau^k}{\tau^k}\delta\ge1$  for all $k\ge0$, and
$\sum_{k=0}^\infty a^k<\infty$, $\sum_{k=0}^\infty b^k<\infty$.
Then, we have:
\begin{itemize}
\item[(a)]  $\lim_{k\to\infty} F(\bar x^k)$
exists and
\[\lim_{k\to\infty}\|\bar y^k\|=
\lim_{k\to\infty}\|x_i^k - \bar x^k\|= \lim_{k\to\infty}\|y_i^k -
\bar y^k\|=0, \: \forall i \]
\item[(b)] If $\{\gamma^k\}$ satisfies $\sum_{k=0}^\infty \gamma^k =\infty$
and  $\lim_{k\rightarrow\infty}\gamma^k>0$, then
$\lim_{k\rightarrow\infty}\|\nabla F(\bar x^k)\|=0$. Moreover, if
$\{\bar x^k\}$ is bounded, then every accumulation point of  $\{\bar
x^k\}$ is an optimal solution, and $\lim_{k\rightarrow\infty}
F(x_i^k)= F(\theta^*)$  for all $i$.
\end{itemize}
\end{Proposition 1}

\begin{proof}
\noindent(a) The  idea is to show that Lemma~\ref{th-dsystem}
applies. Setting up the equation $\pi^T V=\pi^T$, we have $
 (1-\eta)\pi_1+(1-\eta)^2(1-c)(1-\delta) \pi_3=(1-\eta)\pi_2$ and
$\pi_2=(1-\eta)(1-c)\pi_3$. Dividing the first equation with
$1-\eta$, we find
\[\pi_1+(1-\eta)(1-c)(1-\delta) \pi_3=\pi_2
\]
which in view of $\pi_2=(1-\eta)(1-c)\pi_3$ implies
\[
\pi_1 +(1-\delta)\pi_2=\pi_2\quad\implies\quad\pi_1=\delta\pi_2
\]
Thus, for the vector $\pi$  satisfying $\pi^T=V\pi^T$, we have
\be\label{eq-pivec} \pi_1=\delta\pi_2,\qquad\pi_2=(1-\eta)(1-c)\pi_3
\ee Hence, we can find such a vector $\pi$ with $\pi>0$. We next
verify that such a vector also satisfies $\pi^TC>0$. We have
$\pi^TC=[\g^k\pi_1,
\frac{1-\tau^k}{\tau^k}\pi_1-(1-\eta)^2(1-c)\pi_3]$, which,
under~\eqref{eq-pivec},  implies  the second coordinate of $\pi^TC$
satisfying
$[\pi^TC]_2=\frac{1-\tau^k}{\tau^k}\delta\pi_2-(1-\eta)\pi_2
=\left(\frac{1-\tau^k}{\tau^k}\delta-1+\eta\right)\pi_2 $.

 The
condition $\frac{1-\tau^k}{\tau^k}\delta\ge1$  implies
$[\pi^TC]_2\ge\eta\pi_2>0$. Thus, Lemma~\ref{th-dsystem}'s
 conditions
 are satisfied, and it follows that for the
three elements of $\bv^k$, i.e., $\bv_1^k$, $\bv_2^k$, and
$\bv_3^k$, we have that \be\label{eq-limit-exist} \lim_{k\to\infty}
\pi_1\bv^k_1+\pi_2 \bv_2^k+\pi_3\bv_3^k \ee
 exists and
$\sum_{k=0}^\infty\pi^TC\bu^k<\infty$ holds with $\bu^k= [\|\nabla
F(\bar x^k)\|^2,\ \|\bar y^k\|^2]^T$. Since $\pi^TC\geq [\g^k\pi_1,\
\eta\pi_2]$, one has \be\label{eq-sumfinite}
\sum_{k=0}^\infty\g^k\|\nabla F(\bar x^k)\|^2<\infty, \qquad
\sum_{k=0}^\infty\|\bar y^k\|^2<\infty\ee
and hence,
\be\label{eq-to0} \lim_{k\to\infty}\|\bar y^k\|=0\ee If we had that
$\sum_{i=1}^m\|x_i^k-\bar x^k\|^2$ and $\sum_{i=1}^m\|y_i^k-\bar
y^k\|^2$ are convergent, then it would follow
from~\eqref{eq-limit-exist} that the limit $\lim_{k\to\infty} F(\bar
x^k)$ exist.

Now, we focus on proving that both $\bv_2^k=\sum_{i=1}^m\|x_i^k-\bar
x^k\|^2$ and $\bv_3^k=\sum_{i=1}^m\|y_i^k-\bar y^k\|^2$ converge to
0. The  idea is to show that we can apply Lemma~\ref{th-dsystem2}.
By focusing on the elements $\bv_2^k$ and $\bv_3^k$, from
\eqref{eq-fine} we have \bas \left[\begin{array}{c} \bv_2^{k+1}\cr
\bv_3^{k+1}\end{array} \right] \le \left(\tilde V + a^k {\bf 1}{\bf
1}^T\right) \left[\begin{array}{c} \bv_2^k\cr
\bv_3^k\end{array}\right] +\hat b^k{\bf 1} + \left[\begin{array}{c}
 0\cr
\hat{c}^k\end{array}\right] \eas where $\hat b^k = b^k+a^k\nu(F(\bar
x^k)-F(\theta^*))$,  $
 \hat{c}^k =(1-\eta)^2 (1-c)\|\bar y^k\|^2$, and $ \tilde
V=\left[\begin{array}{ccc}
 \eta &\frac{1}{1-\eta}\cr
(1-\eta)^2 (1-c)(1-\delta) & c\end{array}\right]$.

 By separating the first term on the right hand side and bounding
the last vector   by $(1-\eta)^2(1-c)\|\bar y^k\|^2{\bf 1}$, we
obtain \ba\label{eq-finer} \left[\begin{array}{c} \bv_2^{k+1}\cr
\bv_3^{k+1}\end{array} \right] \le \tilde V \left[\begin{array}{c}
\bv_2^k\cr \bv_3^k\end{array}\right] +\tilde b^k{\bf 1} \ea where
\begin{equation}\label{eq:b^k^tilde}
\begin{aligned}
&\tilde b^k= b^k+(1-\eta)^2(1-c)\|\bar y^k\|^2+a^k\times\\
&\left(\nu(F(\bar x^k) - F(\theta^*)) + \sum_{i=1}^m\|x_i^k-\bar
x^k\|^2 +\sum_{i=1}^m\|y_i^k-\bar y^k\|^2\right)
\end{aligned}
\end{equation}

To apply Lemma~\ref{th-dsystem2}, we show that the equation $\pi^T
\tilde V=(1-\a) \pi^T$ has a solution in $\pi=(\pi_2,\pi_3)$ with
$(\pi_2,\pi_3)>0$ and $\a\in(0,1)$. Note that, if we have such a
solution, then we will let $\a^k=\a>0$ for all $k$, so that the
condition $\sum_{k=0}^\infty\a^k=\infty$ of Lemma~\ref{th-dsystem2}
will be satisfied. In this case, the condition
$\lim_{k\to\infty}\tilde b^k/\a^k=0$ of Lemma~\ref{th-dsystem2} will
also be satisfied. This is because by our assumption on the
sequences $\{a^k\}$ and $\{b^k\}$, it follows that
$\lim_{k\to\infty}a^k=0$ and $\lim_{k\to\infty}b^k=0$. We also have
 $\lim_{k\to\infty}\|\bar y^k\|=0$ (see \eqref{eq-to0}). Moreover, in
view of relation~\eqref{eq-limit-exist}, the sequences $\{F(\bar
x^k)-F(\theta^*)\}$, $\sum_{i=1}^m\|x_i^k-\bar x^k\|^2$, and
$\sum_{i=1}^m\|y_i^k-\bar y^k\|^2$ are bounded. Hence, it follows
that $\tilde b^k$ defined in (\ref{eq:b^k^tilde}) will converge to 0
as $k$ tends to infinity. Thus, all the conditions of
Lemma~\ref{th-dsystem2} will be satisfied.

It remains to show that the system of equations $\pi^T \tilde
V=(1-\a) \pi^T$ has a solution in $\pi=(\pi_2,\pi_3)$  with
$(\pi_2,\pi_3)>0$ and $\a\in(0,1)$. The system is equivalent to
\[(1-\eta)^2 (1-c)(1-\delta)\pi_3\hspace{-0.1cm}=\hspace{-0.1cm}(1-\eta-\a)\pi_2,\:
\pi_2\hspace{-0.1cm}=\hspace{-0.1cm}(1-\eta)(1-c-\a)\pi_3\] which
gives $\pi_2>0$ with arbitrary $\pi_3>0$, and imposes that $\a$
satisfies $
 (1-\eta)(1-c)(1-\delta)=(1-\eta-\a)(1-c-\a)$, or equivalently,
\be\label{eq-quadr} \a^2-(2-\eta-c)\alpha+\delta(1-\eta)(1-c)=0\ee
Letting $\psi(\a)=\a^2-(2-\eta-c)\alpha+\delta(1-\eta)(1-c)$ for all
$\a\in\mathbb{R}$, we note that $\psi(\cdot)$ is strongly convex and
its minimum is attained at $ \a_0=\frac{1}{2}(2-\eta-c)$.

 For the
minimum value we have
\[
\begin{aligned}
\psi(\a_0)
&=\textstyle\frac{1}{4} (2-\eta-c)^2 -\frac{1}{2}  (2-\eta-c)^2 +\delta(1-\eta)(1-c)\\
&=-\textstyle\frac{1}{4}  (2-\eta-c)^2 +\delta(1-\eta)(1-c)\\
&=-\textstyle\frac{1}{4}  (1-\eta+1-c)^2  + \delta(1-\eta)(1-c) \\
\end{aligned}
\]

Since $\delta<1$, it follows that $ \psi(\a_0)  <
-\frac{(1-\eta+1-c)^2}{4} +  (1-\eta)(1-c)    =
-\frac{((1-\eta)-(1-c))^2}{4}
 \leq 0
$. We also have $\psi(0)= \delta(1-\eta)(1-c)>0$ since $\delta>0$
and $c,\eta\in(0,1)$. Thus, we have $\psi(0)>0$ and $\psi(\a_0)<0$,
implying that there exists some $\a^*\in(0,\a_0)$ satisfying
$\psi(\a^*)=0$ with $\a_0=\frac{2-\eta-c}{2}$. Since
$c,\eta\in(0,1)$, we have $\a_0\in(0,1)$. Hence,  \eqref{eq-quadr}
has a solution $\a^*\in(0,1)$. So there is a vector $\pi>0$ and
$\a\in(0,1)$ that satisfy $\pi^T \tilde V=(1-\a) \pi^T$, and we can
apply Lemma~\ref{th-dsystem2} with $\a_k=\a$ for all $k$. By
Lemma~\ref{th-dsystem2}, we have $\lim_{k\to\infty}\|x_i^k - \bar
x^k\|=0$ and $\lim_{k\to\infty}\|y_i^k - \bar y^k\|=0$.

\noindent(b) Since $\sum_{k=0}^\infty\g^k\|\nabla F(\bar
x^k)\|^2<\infty$ (see~\eqref{eq-sumfinite}), from $\sum_{k=0}^\infty
\gamma^k =\infty$ and   $\lim_{k\to\infty}\g^k>0$, it follows
$\lim_{k\to\infty}\|\nabla F(\bar x^k)\|=0$.

Now, if $\{\bar x^k\}$ is bounded, then it has accumulation points.
Let $\{\bar x^{k_i}\}$ be a sub-sequence such that
$\lim_{i\to\infty}\|\nabla F(\bar x^{k_i})\|=0$. Without loss of
generality, we may assume that $\{\bar x^{k_i}\}$ is convergent, for
otherwise we would choose a sub-sequence of $\{\bar x^{k_i}\}$. Let
$\lim_{i\to\infty}\bar x^{k_i}=\hat x$. Then, by continuity of the
gradient $\nabla F(\cdot)$, it follows $\nabla F(\hat x)=0$,
implying that $\hat x$ is an optimal point. Since $F$ is continuous,
it follows   $\lim_{i\to\infty} F(\bar x^{k_i}) = F(\hat
x)=F(\theta^*)$. By part (a),  $\lim_{k\to\infty} F(\bar x^k)$
exists, so we must have $\lim_{k\to\infty} F(\bar x^k)=F(\theta^*)$.

Finally,  by part (a) we have $\lim_{k\to\infty}\|x_i^k-\bar
x^k\|^2=0$ for every $i$. Thus, it follows that each  $\{x_i^k\}$
has the same accumulation points as $\{\bar x^k\}$, implying by
continuity of the objective function $F$ that $\lim_{k\to\infty}
F(x^k_i)=F(\theta^*)$ for all $i$.
\end{proof}

\begin{Theorem 1}\label{theorem:PDGD_NDS_convergence}
Under  Assumption 1 and Assumption 2, if there exists some $T\geq 0$ such that
for all $k\geq T$, the stepsize vector  $ \lambda^k = [\lambda_1^k,\cdots,\lambda_m^k]^T $ ( with all elements non-negative)  satisfies
\[\sum_{k=T}^\infty \bar\lambda^k =\infty,\:
\sum_{k=T}^\infty \|\lambda^{k+1}-\lambda^k\|^2 <\infty,\:
\sum_{k=T}^\infty\frac{\|\lambda^k-\bar\lambda^k{\bf 1}\|^2}
{\bar\lambda^k}<\infty\] with $\bar\lambda^k=\frac{\sum_{i=1}^{m}\lambda_i^k}{m}$, and
\[
\begin{aligned}
&\frac{2L}{m\bar\lambda^k}(\lambda_{\max}^k)^2\le 1-\eta,
\bar\lambda^k\le \frac{\delta}{1+\delta},    \eta
+\frac{6m^2L^2}{1-\eta}\|\lambda^{k+1}\|^2 \le c,\\
& \max\{m^3r^2,m^2\} 6L^2\|\lambda^{k+1}\|^2\le
(1-\eta)^3(1-c)(1-\delta)
\end{aligned}
\]
for some $\delta\in(0,1)$,
$c\in(0,1)$,
 then, the results of Proposition~\ref{th-main2}
hold for the proposed PDG-NDS.
\end{Theorem 1}
\begin{proof}
The  idea is to prove that we can establish the relationship in
(\ref{eq-fine}). To this end, we divide the derivation into four
steps: in Step I, Step II, and Step III, we establish  relationships
for $\frac{2}{L}\left(F(\bar x^{k+1}) - F(\theta^*)\right)$,
$\sum_{i=1}^m \|x_i^{k+1} -\bar x^{k+1}\|^2$, and $\sum_{i=1}^m
\|y_i^{k+1}-\bar y^{k+1}\|^2$, respectively, and in Step IV, we
prove that (\ref{eq-fine}) holds. { To help the exposition of the main idea, we put Step I, Step II, and Step III in Appendix B, and only give the  derivation of Step IV here. } 

Step IV: We summarize the relationships obtained in Steps I-III in Appendix B and
prove the theorem. Defining $\bv^k=\big[\frac{2}{L}(F(\bar
x^{k+1})-F(\theta^*)),\sum_{i=1}^m\|x_i^{k+1}-\bar
x^{k+1}\|^2,\sum_{i=1}^m\|y_i^{k+1}-\bar y^{k+1}\|^2\big]^T$, we
have the following relations from (\ref{eq:StepI_final}),
(\ref{eq:StepII_final_b}), and (\ref{eq:stepIII_final}) in Appendix B:
\begin{equation}\label{eq:stacked}
\begin{aligned}
\bv^{k+1}\leq(V^k+A^k)\bv^k-C^k\left[\begin{array}{c} \left\|\nabla
F(\bar x^k)\right\|^2\cr \|\bar y^k\|^2
\end{array}\right]+B^k
\end{aligned}
\end{equation}
where
\[
\begin{aligned}
&V^k=\left[\begin{array}{ccc} 1 &
\frac{2L}{m\bar\lambda^k}(\lambda_{\max}^k)^2& 0\cr 0& \eta &
\frac{1}{1-\eta}\cr 0&
\frac{6m^2r^2L^2}{1-\eta}\left\|\lambda^{k+1}\right\|^2&
\eta+\frac{6m^2L^2}{1-\eta}\|\lambda^{k+1}\|^2
\end{array}\right],\\
&A^k=\left[\begin{array}{ccc} \frac{c_1}{\bar
\lambda^k}\|\lambda^k-\bar\lambda^k{\bf 1}\|^2& 0 &0\cr 0& 0 &0\cr
\frac{8m^3L^2}{1-\eta}\left\|\lambda^{k+1}-\lambda^k\right\|^2  &
\frac{4m^2L^2}{1-\eta}\left\|\lambda^{k+1}-\lambda^k\right\|^2 &0
\end{array}\right],\\
&C^k=\left[\begin{array}{cc}  \frac{\bar\lambda^k}{L} &
 \frac{1-\bar \lambda^k L}{\bar\lambda^kL}\cr 0&0\cr
0&-\frac{6m^3r^2L^2}{1-\eta}\left\|\lambda^{k+1}\right\|^2
\end{array}\right],\\
&B^k=\left[\begin{array}{c} \frac{c_1}{\bar
\lambda^k}\|\lambda^k-\bar\lambda^k{\bf 1}\|^2 \sum_{i=1}^m\|\nabla
f_i(\theta^*)\|^2\cr 0\cr
\frac{8m^2}{1-\eta}\left\|\lambda^{k+1}-\lambda^k\right\|^2\sum_{i=1}^m\|\nabla
f_i(\theta^*)\|^2
\end{array}\right]
\end{aligned}
\]
Now using the conditions of the theorem, we bound  the entries in
$V^k$, $A^k$, $C^k$, and $B^k$. It can be seen that
\begin{equation}\label{eq:A^k}
A^k\leq a^k{\bf1}{\bf 1}^T,\quad B^k\leq b^k{\bf 1}
\end{equation}
hold where \be\label{def-ak} a^k=\max\left\{\frac{c_1}{\bar
\lambda^k}\|\lambda^k-\bar\lambda^k{\bf 1}\|^2,
\frac{8m^3L^2}{1-\eta}\left\|\lambda^{k+1}-\lambda^k\right\|^2\right\}\ee
\begin{equation}
\begin{aligned}
\label{def-bk} b^k=&\max\left\{\frac{c_1}{\bar
\lambda^k}\|\lambda^k-\bar\lambda^k{\bf 1}\|^2,
\frac{8m^2}{1-\eta}\left\|\lambda^{k+1}-\lambda^k\right\|^2\right\}\\
&\times \sum_{i=1}^m\|\nabla f_i(\theta^*)\|^2
\end{aligned}
\end{equation}

Using  $\eta+\frac{6m^2L^2}{1-\eta}\|\lambda^{k+1}\|^2\le c$ with
$c\in(0,1)$, $\frac{2L}{m\bar\lambda^k}(\lambda_{\max}^k)^2\le
1-\eta$, and $m^2r^26L^2\|\lambda^{k+1}\|^2\leq
m^3r^26L^2\|\lambda^{k+1}\|^2\le (1-\eta)^3(1-c)(1-\delta)$ from the
theorem conditions, we can bound $V^k$:
\begin{equation}\label{eq:V^k}
V^k\leq V\triangleq \left[\begin{array}{ccc} 1 & 1-\eta& 0\cr 0&
\eta & \frac{1}{1-\eta}\cr 0& (1-\eta)^2(1-c)(1-\delta)& c
\end{array}\right]
\end{equation}
Furthermore, we can bound $C^k$ using the condition
$\max\{m^3r^2,m^2\} 6L^2\|\lambda^{k+1}\|^2\le
(1-\eta)^3(1-c)(1-\delta)$ which implies
$\frac{m^3r^26L^2}{1-\eta}\|\lambda^{k+1}\|^2\le (1-\eta)^2(1-c)$:
\begin{equation}\label{eq:C^k}
C^k\geq C\triangleq \left[\begin{array}{cc} \frac{\bar\lambda^k}{L}
& \frac{1-\bar \lambda^k L}{\bar\lambda^kL}\cr 0&0\cr
0&-(1-\eta)^2(1-c)
\end{array}\right]
\end{equation}

Combining (\ref{eq:stacked}), (\ref{eq:A^k}), (\ref{eq:V^k}), and
(\ref{eq:C^k}) leads to
\begin{equation}\label{eq:stacked_final}
\begin{aligned}
\bv^{k+1}=(V+a^k{\bf1}{\bf 1}^T)\bv^k-C\left[\begin{array}{c}
\left\|\nabla F(\bar x^k)\right\|^2\cr \|\bar y^k\|^2
\end{array}\right]+b^k{\bf1}
\end{aligned}
\end{equation}

We note that Proposition \ref{th-main2} applies to the  relations in
(\ref{eq:stacked_final}) for $k\ge T$, with
$\g^k=\frac{\bar\lambda^k}{L}$, $\tau^k=\bar \lambda^k L$, and $a^k$
and $b^k$ given by~\eqref{def-ak} and~\eqref{def-bk}, respectively.
By our assumption  $\sum_{k=T}^\infty \|\lambda^{k+1}-\lambda^k\|^2
<\infty$ and $\sum_{k=T}^\infty\frac{\|\lambda^k-\bar\lambda^k{\bf
1}\|^2} {\bar\lambda^k}<\infty$, it follows that $\{a^k\}$ and
$\{b^k\}$ are nonegative and summable. The condition $\bar \lambda^k
L\le \delta/(1+\delta)$ is equivalent to $\bar \lambda^kL +\delta
\bar \lambda^kL\le \delta$, implying  $1\le \frac{\delta(1-\bar
\lambda^kL)}{\bar \lambda^k L}$. Thus, with $\tau^k=\bar \lambda^k
L$, we see that the condition $\frac{1-\tau^k}{\tau^k }\delta\ge 1$
of Proposition~\ref{th-main2} is satisfied for all $k\ge T$.
Additionally, by our assumption $\sum_{k=T}^\infty \bar\lambda^k
=\infty$ we see that the condition of Proposition~\ref{th-main2}(b) also
holds for $k\ge T$. Since the results of Proposition~\ref{th-main2} are
asymptotic, the results remain valid when the starting index is
shifted from $k=0$ to $k=T$, for an arbitrary $T\ge 0$.
\end{proof}

{

\begin{Remark 1}
To enable privacy protection, we let each agent $i$ decide its stepsize $\lambda_i^k$.
   To make sure that the selected stepsizes  satisfy the convergence condition  in the statement of Theorem \ref{theorem:PDGD_NDS_convergence}, all agents  need  to know an upper bound on the Lipschitz constant $L$. It is   worth noting that such  information about $L$ is needed in all decentralized optimization algorithms with a constant stepsize (see, e.g., \cite{shi2015extra,xu2015augmented,xin2018linear,pu2020push}). In fact, in implementations, to satisfy the
condition in   the statement of Theorem \ref{theorem:PDGD_NDS_convergence}, all agents can be given the same baseline value of stepsize $ \lambda $. Then, every agent can set its stepsize $\lambda_i^k$ by deviating from the baseline value in a finite
number of iterations. The indices of these iterations are private
to individual agents. As long as the deviation in each of these
iterations is finite, the heterogeneity condition in   the statement of Theorem \ref{theorem:PDGD_NDS_convergence}  will
be satisfied.
\end{Remark 1}
}


\begin{Remark 1}\label{remark:conversion}
It is worth noting that although some gradient-tracking based
algorithms can be reduced to the $x$-variable only form by
eliminating the auxiliary variable, and hence share one variable in
interaction, such a conversion is infeasible when the stepsizes are
heterogeneous and not shared across agents (for the purpose of,
e.g., privacy preservation). For example, the Aug-DGM algorithm in
\cite{xu2015augmented} has the following form
\[
\left\{
\begin{aligned}
x^{k+1}&=W(x^k-\Lambda y^k)\\
y^{k+1}&=W(y^k+g^{k+1}-g^k)
\end{aligned}
\right.
\]
Although we can eliminate the $y$ variable and convert it to
\[
\begin{aligned}
x^{k+2}&=(W+W\Lambda W\Lambda^{-1}W^{-1})x^{k+1}-W\Lambda
W\Lambda^{-1}x^k\\
&\qquad -W\Lambda W(g^{k+1}-g^k)
\end{aligned}
\]
we cannot   let agent $j$ share $(W+W\Lambda
W\Lambda^{-1}W^{-1})_{ij}x_j^{k+1}-(W\Lambda
W\Lambda^{-1})_{ij}x_j^k-(W\Lambda W)_{ij}(g_j^{k+1}-g_j^k)$ in each
iteration when the stepsizes are   not shared across agents (for the
purpose of, e.g., privacy preservation). This is because calculating
$(W+W\Lambda W\Lambda^{-1}W^{-1})_{ij}$ and $(W\Lambda
W\Lambda^{-1})_{ij}$ requires agent $j$ to know all stepsizes
$\Lambda$, which however, were assumed to be private to individual
agents. Therefore, even though some existing gradient-tracking based
algorithms can use heterogeneous stepsizes to hide information, they
have to exchange two messages between interacting agents. In fact,
privacy  enabled in this way is   quite weak, as detailed in Sec.
\ref{sec:privacy_analysis_NDS}.
\end{Remark 1}
\subsection{Privacy analysis}\label{sec:privacy_analysis_NDS}

Next we show that  PDG-NDS  is able to provide every participating
agent the privacy defined in Definition \ref{definition:privacy}
against both honest-but-curious and eavesdropping adversaries.

\begin{Theorem 3}\label{theo:privacy_diminishing}
In the presence of honest-but-curious or eavesdropping adversaries,
PDG-NDS can protect the privacy of all participating agents defined
in Definition \ref{definition:privacy}.
\end{Theorem 3}

\begin{proof}

{ Without loss of generality, we first consider the protection of the gradient of agent $i$ at any single time instant $k$, and then show that the argument also applies to any finite number of time instants (iterations).}  When the gradient  is $g^k_i$, we
represent the information that agent $i$ shares with others in
PDG-DS as $\mathcal{I}_i$. According to
Definition \ref{definition:privacy}, we  have to prove that at some
iteration $k$, if the gradient is altered from $g^k_i$ to $\hat{g}^k_i=e^{\zeta^k}g^k_i$ with $\zeta^k$ difference from $g^k_i$ according to the log-scale metric in (\ref{eq:difference}), the corresponding shared
information $\hat{\mathcal{I}}_i$ of agent $i$ when running the
algorithm  could be identical to $\mathcal{I}_i$ under any
$\zeta^k>0$.

According to Algorithm PDG-NDS, agent $i$ shares the following
information in decentralized optimization:
\[
\mathcal{I}_i=\mathcal{I}_i^{\rm sent}\bigcup\mathcal{I}_i^{\rm
public}
\]
with $ \mathcal{I}_i^{\rm sent}
=\left\{v^k_{ji}|k=1,2,\cdots\right\} $,
 $
v^k_{ji}
=2w_{ji}x_i^{k-1}+\{W^2\}_{ji}x_i^{k-2}-b_{ji}^{k-2}(\lambda_i^{k-1}g_i^{k-1}-\lambda_i^{k-2}g_i^{k-2})
$ , and  $\mathcal{I}_i^{\rm public}=\left\{W\bigcup
\sum_{j\in\mathbb{N}_i}b^k_{ji}=1|k=0,1,\cdots\right\} $.

 It can be
obtained that when the gradient is altered to
$\hat{g}^k_i=e^{\zeta^k}g^k_i$, the difference defined in
(\ref{eq:difference}) is $\zeta^k$. However, in this case, if we set
the stepsize $\hat{\lambda}_i^k$ to
$\hat{\lambda}_i^k=e^{-\zeta^k}\lambda_i^k$, then the corresponding
shared information will still be $v_{ji}^k$. Since other parameters
are not changed and changing the stepsize {  from $\lambda_i^k$} to
$\hat{\lambda}_i^k=e^{-\zeta^k}\lambda_i^k$ at $k$ will not violate
the summable stepsize heterogeneity conditions in Theorem
\ref{theorem:PDGD_NDS_convergence} for any given $\zeta^k<\infty$,
  convergence to
the optimal solution will still be guaranteed. {  Similarly, if the gradient of agent $i$ is altered at iteration  $k$ and iteration $k+1$ to $\hat{g}^k_i=e^{\zeta^k}g^k_i$ and $\hat{g}^{k+1}_i=e^{\zeta^{k+1}}g^{k+1}_i$, respectively, these alterations can be covered by a stepsize alteration of $\hat{\lambda}_i^k=e^{-\zeta^k}\lambda_i^k$ at iteration $k$ and $\hat{\lambda}_i^{k+1}=e^{-\zeta^{k+1}}\lambda_i^{k+1}$ at iteration $k+1$.  Given that the
convergence condition in the statement of Theorem  \ref{theorem:PDGD_NDS_convergence}
  allows the stepsize of agent $i$ to change by any finite amount for any finite number of iterations, one can obtain that the variations of every agent's gradients in any number of iterations can be completely covered by the flexibility in changing the agent's stepsizes in these iterations, as long as the  number of these iterations is finite.}  
Therefore,
privacy of the gradient information of any agent will be protected
when running PDG-NDS. { It is worth noting that the perturbation does not violate the convexity and Lipschitz conditions in Assumption 2. This is because in order for an adversary to check if Assumption 2 is violated, it has to know $x_i^{k}$, which, however, is not available to adversaries: before convergence, $x_i^k$ is inaccessible to the adversary because the information shared by agent $i$ is $ 2w_{ji}x_i^{k-1}+\{W^2\}_{ji}x_i^{k-2} -b_{ji}^{k-2}(\lambda_i^{k-1}g_i^{k-1}-\lambda_i^{k-2}g_i^{k-2})$, avoiding $x_i^{k-1}$ and $x_i^{k-2}$ from being inferrable.  In fact, even with the gradient $g_i^{k-1}$ and $g_i^{k-2}$ unchanged, the value of observation  $ 2w_{ji}x_i^{k-1}+\{W^2\}_{ji}x_i^{k-2} -b_{ji}^{k-2}(\lambda_i^{k-1}g_i^{k-1}-\lambda_i^{k-2}g_i^{k-2})$ can be changed by an arbitrary finite value by changing the stepsize $\lambda_i^{k-1}$ or $\lambda_i^{k-2}$.
 Therefore,  before convergence, an adversary cannot use Assumption 2 to confine the change in observed values and further confine the change in the gradient. After convergence, the perturbation does not violate the convexity and Lipschitz conditions, either.  In fact, although $x_i^k$ becomes accessible to the adversary after convergence,   gradient information  is eliminated in adversary's observation (the shared information $ 2w_{ji}x_i^{k-1}+\{W^2\}_{ji}x_i^{k-2} -b_{ji}^{k-2}(\lambda_i^{k-1}g_i^{k-1}-\lambda_i^{k-2}g_i^{k-2})$ becomes  $ 2w_{ji}x_i^{k-1}+\{W^2\}_{ji}x_i^{k-2}$ because $\lambda_i^{k-1}$ and $\lambda_i^{k-2}$ converge  to the same constant value). So after convergence,   the adversary still cannot use Assumption 2 to confine  changes in gradients.}
\end{proof}

{
\begin{Remark 1}\label{re:proof_convergence_algorithm2}
Even after convergence when $g_i^k$ becomes a constant, an adversary still cannot infer gradients from shared messages in PDG-NDS. More specifically, when $g_i^k$ converges to a constant value, the stepsize $\lambda_i^k$  also converges to a constant value, which completely eliminates the information of $g_i^k$  in observed information (the observed information becomes $2w_{ji}x_i^{k-1}+\{W^2\}_{ji}x_i^{k-2}$ after convergence). This can also be understood intuitively as follows: Even if the adversary can collect  $T\rightarrow\infty$  observations $ 2w_{ji}x_i^{k-1}+\{W^2\}_{ji}x_i^{k-2} -b_{ji}^{k-2}(\lambda_i^{k-1}g_i^{k-1}-\lambda_i^{k-2}g_i^{k-2})$ in the neighborhood of the optimal point and establish  a system  of $T$ equations to solve for $g_i^{k-1}$ and $g_i^{k-2}$ (which can be   viewed to be approximately time-invariant  in the neighborhood of the optimal point), the number of unknowns $b_{ji}^{k-2}$, $\lambda_i^{k-1}$,  $\lambda_i^{k-2}$, $g_i^{k-1}$ and $g_i^{k-2}$ in the system  of $T$ equations is $5T$ (even if we view $\lambda_i^{k-1}$ and $\lambda_i^{k-2}$ to be approximately constant and equal to each other,  $g_i^{k-1}$ and $g_i^{k-2}$ to be constant and equal to each other, in the neighborhood of the optimal point, the number of unknowns is still $T+2$), which makes it impossible for the adversary to solve for $g_i^{k-1}$ or $g_i^{k-2}$ using the system  of $T$ equations established from $T$  observations.
\end{Remark 1}

}
\begin{Remark 1}
Similar to PDG-DS,  our PDG-NDS algorithm can protect the privacy of
every participating agent against honest-but-curious and
eavesdropping adversaries without any constraint on the interaction
topology.
\end{Remark 1}

\begin{Remark 1}
In line with the discussions in the previous remark, PDG-NDS can
protect the defined privacy of participating agents irrespective of
the number of adversaries and existence of collusion among
adversaries.
\end{Remark 1}

 We next show that {  directly making stepsize and coupling matrices time-varying in existing gradient-tracking based
 decentralized optimization algorithms cannot   provide
 the defined privacy. We use  the AB algorithm proposed in \cite{xin2018linear} as an example
to show this since it allows  a column-stochastic coupling matrix, which allows individual agents to keep their coupling coefficients private. The AB algorithm has the following form \cite{xin2018linear}:
 \[
\left\{
\begin{aligned}
x^{k+1}&= Rx_j^k-\lambda y^k\\
y^{k+1}&= C(y^k+g^{k+1}-g^k)
\end{aligned}
\right.
\]
where   $R=\{r_{ij}\}$ is row-stochastic and
$C=\{c_{ij}\}$ is column-stochastic.

Directly making its stepsize and coupling coefficients time-varying leads to the following algorithm (we also introduce heterogeneity in the stepsize):}
 \[
\left\{
\begin{aligned}
x^{k+1}&= R^kx_j^k-\Lambda^k y^k\\
y^{k+1}&= C^k(y^k+g^{k+1}-g^k)
\end{aligned}
\right.
\]
where  $R^k=\{r^k_{ij}\}$   should be row-stochastic and
$C^k=\{c^k_{ij}\}$ should be column-stochastic. At each iteration
$k$, an agent $j$ will share $x_j^k$ and
$c^k_{ij}(y_j^k+g_j^{k+1}-g^k_j)$ with its neighboring agent $i$.
Also, all agents initialize as $y_i^0=g_i^0$.

Because for all $i\in\mathbb{N}_{j}$,  $c^k_{ij}$ are generated by
agent $j$, it seems that agent $j$ can keep $c^k_{ij}$ confidential
and hence uses them to cover  shared information
$c^k_{ij}(y_j^k+g_j^{k+1}-g_j^k)$. Next, we show that this is not
true.

We consider the case where agent $i$ is the only neighbor of agent
$j$. In this case, agent $i$ knows agent $j$'s update rule
\[
y_j^{k+1}=
c^k_{jj}(y_j^k+g_j^{k+1}-g_j^k)+c^k_{ji}(y_i^k+g_i^{k+1}-g_i^k)
\]
Using the fact $c^k_{jj}+c^k_{ij}=1$, the above update rule can be
rewritten as
\[
y_j^{k+1}=
(1-c^k_{ij})(y_j^k+g_j^{k+1}-g_j^k)+c^k_{ji}(y_i^k+g_i^{k+1}-g_i^k)
\]
or
\begin{equation}\label{eq:privacy_analysis_AB}
\begin{aligned}
&y_j^{k+1}-y_j^k=g_j^{k+1}-g_j^k\\
&\qquad
 -c^k_{ij} (y_j^k+g_j^{k+1}-g_j^k)+c^k_{ji}(y_i^k+g_i^{k+1}-g_i^k)
 \end{aligned}
\end{equation}
Note that $c^k_{ij} (y_j^k+g_j^{k+1}-g_j^k)$ is shared with agent
$i$ by agent $j$ and $c^k_{ji}(y_i^k+g_i^{k+1}-g_i^k)$ is generated
by agent $i$, hence the last two terms on the right hand side of
(\ref{eq:privacy_analysis_AB}) are known to agent $i$. We represent
$-c^k_{ij} (y_j^k+g_j^{k+1}-g_j^k)+c^k_{ji}(y_i^k+g_i^{k+1}-g_i^k)$
as $m_i^k$ and sum (\ref{eq:privacy_analysis_AB}) from $k=0$ to $t$
to obtain
$
y_j^{t+1}=g_{j}^{t+1}+\sum_{k=0}^{t}m_i^k
$
where we used the relationship $y_j^0=g_j^0$.

When $t\rightarrow\infty$, we have $y_j^t\rightarrow 0$ and
$x_i^k\rightarrow x_j^k$ in the AB algorithm, resulting in
$g_{j}^{t+1}=-\sum_{k=0}^{t}m_i^k$. Therefore, agent $i$ can  infer
the gradient of agent $j$ based on its accessible information
 $m_i^k$. Note that  the above derivation is independent of
the evolution of $x^k$ and stepsize $\Lambda^k$, so the same privacy
leakage will occur even if the stepzies are uncoordinated (not
shared across agents) such as in
\cite{xu2015augmented,nedic2017geometrically}.

One may wonder if we can reduce the AB algorithm to the $x$-variable
only form to avoid information leakage. Given that when the
stepsizes are heterogeneous and not shared across agents, such
reduction is impossible, as detailed in Remark
\ref{remark:conversion}, we only consider the homogeneous stepsize
case. In fact, after eliminating the $y^k$ variable, the AB
algorithm reduces to
\begin{equation}\label{eq:AB_x_only}
x^{k+2}=(R^k+C^k)x^{k+1}-C^kR^kx^k-\lambda C^k(g^{k+1}-g^k)
\end{equation}
Note that for privacy-preserving purposes, agent $j$ should keep
$c_{ij}^k$ private and send
$(R^k+C^k)_{ij}x_j^{k+1}-(C^kR^k)_{ij}x_j^k-\lambda
C^k(g_j^{k+1}-g_j^{k})$ to agent $i$ where $(\cdot)_{ij}$ represents
the $(i,j)$th element of a matrix. Given $(C^kR^k)_{ij}=\sum_{p=1}^m
c^k_{ip}r^k_{pj}$, agent $j$ has to know all elements of  $C^k$ to
implement the algorithm in $x$-variable only form
(\ref{eq:AB_x_only}), which contradicts the assumption that the
elements of $C^k$ are kept private  to cover information. In other
words, if the column-stochastic matrix $C^k$ is used to cover
information, the AB algorithm cannot be implemented in an
$x$-variable only form. In summary, gradient-tracking based
decentralized optimization algorithms cannot be used to enable the
privacy defined in this paper even under time-varying coupling
weights and heterogeneous stepsizes.

\section{Numerical Simulations}
We use numerical simulations to illustrate the effectiveness of the
proposed algorithms. We consider the canonical distributed
estimation problem where a sensor network of $m$ sensors are used to
collectively estimate an unknown parameter $\theta\in\mathbb{R}^d$.
Each sensor has a noisy measurement of the parameter
$z_i=M_i\theta+w_i$ where $M_i\in\mathbb{R}^{s\times d}$ is the
measurement matrix and $w_i$ is   Gaussian noise. The maximum
likelihood estimation problem can be formulated as the decentralized
optimization problem (\ref{eq:optimization_formulation1}) with each
$f_i$ given by $
f_i(\theta)=\|z_i-M_i\theta\|^2+\sigma_i\|\theta\|^2 $ where
$\sigma_i\geq 0$ is the regularization parameter
\cite{xu2017convergence}.

\begin{figure}
    \begin{center}
        \includegraphics[width=0.25\textwidth]{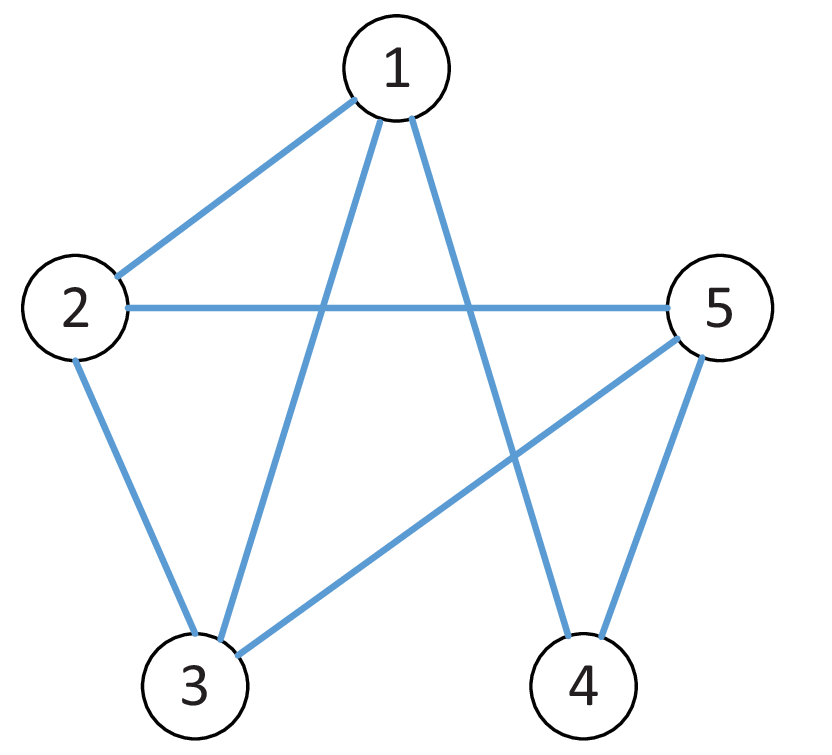}
    \end{center}
    \caption{The interaction topology of the network.}
    \label{fig:topology}
\end{figure}

\begin{figure}
    \begin{center}
        \includegraphics[width=0.5\textwidth]{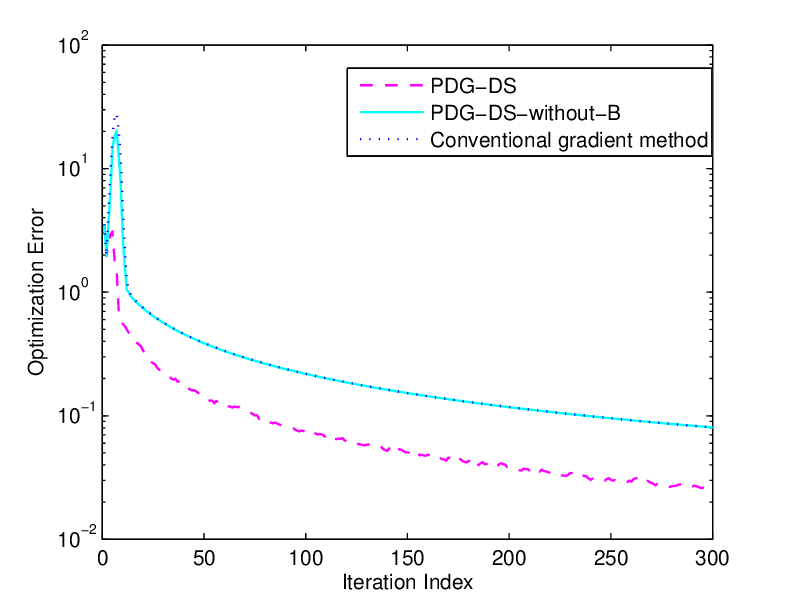}
    \end{center}
    \caption{Comparison of  PDG-DS   with the conventional decentralized gradient method under diminishing-stepsizes.}
    \label{fig:dinimishing_stepsize}
\end{figure}

We considered a network of $m=5$ sensors interacting on the graph
depicted in Fig. \ref{fig:topology}. We set $s=3$ and $d=2$. To
evaluate the performance of our PDG-DS algorithm, we set the
stepsize of agent $i$ as $\lambda_i^k=\frac{1-\varrho_i^k/k^2}{k}$
where $\varrho_i^k$ was randomly chosen by agent $i$ from the
interval $[0,\,1]$ for each iteration. Given that different agents
$i$ chose $\varrho_i^k$ independently, the stepsizes are
heterogeneous across the agents. Each agent $i$ also chose
$b^k_{ji}$ for all $ j\in \mathbb{N}_i$ randomly and independently
of each other under the sum-one condition (to make $B^k$
column-stochastic). The evolution of the optimization error of
PDG-DS is given by the dashed magenta line in Fig.
\ref{fig:dinimishing_stepsize}. When we fixed the $B^k$ matrix to an
identity matrix, we also obtained convergence to the optimal
solution, which is illustrated by the solid green line in Fig.
\ref{fig:dinimishing_stepsize}. The evolution of the conventional
decentralized gradient  algorithm
(\ref{eq:conventional_gradient_descent}) under homogeneous
diminishing stepsize $\frac{1}{k}$ is also presented in Fig.
\ref{fig:dinimishing_stepsize} by the dotted blue line for
comparison. It can be seen that our PDG-DS algorithm has a
comparable (in fact faster) convergence speed with the conventional
decentralized gradient algorithm which does not take privacy into
account. Furthermore, comparing the case with $B^k$ and the case
without $B^k$, we can see that by using the mixing matrix $B^k$ in
PDG-DS, we get faster convergence. This is intuitive since the $B^k$
matrix enhances mixing of information across the agents.

\begin{figure}
\vspace{-0.3cm}
    \begin{center}
        \includegraphics[width=0.5\textwidth]{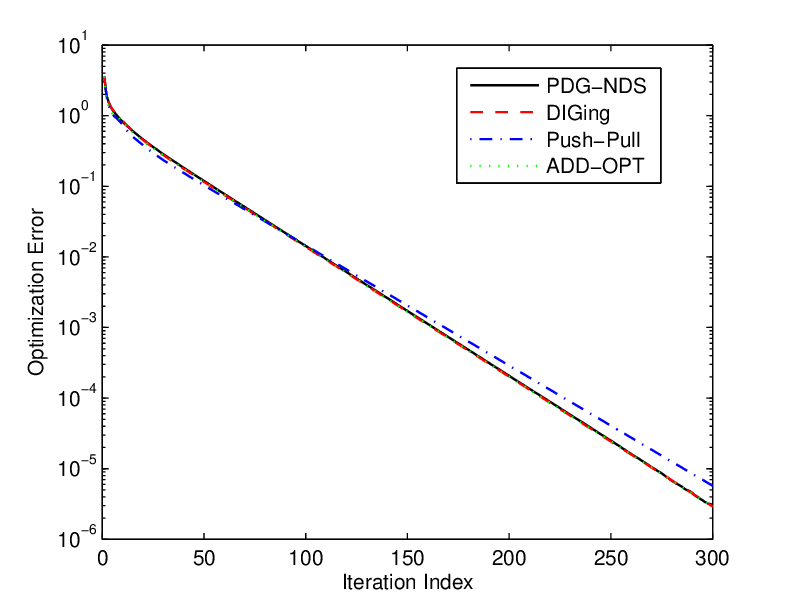}
    \end{center}
    \caption{Comparison of  PDG-NDS  with some gradient-tracking based decentralized optimization algorithms under non-diminishing-stepsizes.}
    \label{fig:non-diminishing_stepsize}
\end{figure}

We also simulated our PDG-NDS algorithm  with non-diminishing
stepsizes. More specifically, we set the stepsize of agent $i$ to
{ $0.02(1-\frac{\varrho^k_i}{k^2})$}  with $\varrho_i^k$ randomly chosen
by agent $i$ from   $[0,\,1]$ for each iteration. Again, since each
agent $i$ chose $\varrho^k_i$ randomly and independently from each
other, the stepsizes are heterogeneous across the network. Each
agent $i$ also randomly chose $b_{ji}^k$ for all $ j\in
\mathbb{N}_i$
  under the sum-one condition (to make $B^k$ column-stochastic). The evolution of the
optimization error of PDG-NDS is illustrated by the solid line in
Fig. \ref{fig:non-diminishing_stepsize}. For the purpose of
comparison, we also plotted the optimization results of
 gradient-tracking algorithms DIGing \cite{nedic2017achieving},
 Push-Pull \cite{pu2020push}, and ADD-OPT \cite{xi2017add} under
 homogeneous stepsize $0.02$. It can be seen that PDG-NDS provides
 similar convergence performance besides enabling privacy protection.


\section{Conclusions}

This paper proposes two inherently privacy-preserving decentralized
optimization algorithms which can guarantee the privacy of all
participating agents without compromising optimization accuracy.
This is in distinct difference from differential privacy based
approaches which trade optimization accuracy for privacy. The two
algorithms are also efficient in communication and computation in
that they are encryption-free and only require an agent to share one
message with a neighboring agent in every iteration, both in the
diminishing stepsize case (the first algorithm) and the
non-diminishing stepsize case (the second algorithm). Note that
commonly used gradient-tracking based decentralized optimization
algorithms require an agent to share two variables, i.e., both the
optimization variable and the gradient-tracking variable. The two
approaches can protect the privacy of every agent even if all
information  shared by an agent is accessible to an adversary, in
which case most existing accuracy-maintaining privacy-preserving
decentralized optimization solutions fail to provide privacy
protection. In fact, even without considering privacy, the
convergence analyses of the two  algorithms under time-varying
uncoordinated stepsizes
   are also of interest by themselves since existing results only
   consider
   constant or fixed
 heterogeneity in stepsizes.
 Numerical simulation results show that both
approaches have similar convergence speeds compared with their
respective privacy-violating counterparts.
\section*{Appendix A}
\begin{Lemma 1}\label{lem-polyak}(\cite{polyak87}, Lemma 11, page 50)
Let $\{v^k\}$, $\{\alpha^k\}$, and $\{\b^k\}$ be sequences of
nonnegative scalars such that $\sum_{k=0}^\infty \a^k<\infty$,
$\sum_{k=0}^\infty \b^k<\infty$, and $v^{k+1}\le(1+\a^k) v^k +\b^k$
holds  for all $k\ge0$. Then, the sequence $\{v^k\}$ is convergent,
i.e., $\lim_{k\to\infty} v^k=v$ for some $v\ge0$.
\end{Lemma 1}

\begin{Lemma 1}\label{lem-polyak2}(\cite{polyak87}, Lemma 10, page 49)
Let $\{v^k\}$, $\{\a^k\}$, and $\{\b^k\}$ be sequences of
nonnegative scalars such that $\sum_{k=0}^\infty \a^k=\infty$,
$\lim_{k\to\infty}\b^k/\a^k=0$, and  $v^{k+1}\le(1-\a^k) v^k +\b^k$
and $\a^k\le 1$ hold for all $k$. Then, the sequence $\{v^k\}$
converges to 0, i.e., $\lim_{k\to\infty} v^k=0$.
\end{Lemma 1}

\begin{Lemma 1}\label{lem-opt}\cite{nedic2014distributed}
Consider a minimization problem $\min_{z \in \R^d} \phi(z)$, where
$\phi:\mathbb{R}^d\to\mathbb{R}$ is a continuous function. Assume
that the optimal solution set $Z^*$ of the problem is nonempty. Let
$\{z^k\}$ be a sequence such that for any optimal solution $z^*\in
Z^*$ and for all $k\ge0,$
\[\|z^{k+1}-z^*\|^2 \le (1+\a^k)\|z^k - z^*\|^2 - \g^k\left (\phi(z^k) - \phi(z^*)\right) +\b^k,\]
where $\a_k\ge0,$ $\b_k\ge0$, and $\g_k\ge0$ for all $k\ge0$, with
$\sum_{k=0}^\infty \a^k<\infty$, $\sum_{k=0}^\infty \g^k=\infty$,
and $\sum_{k=0}^\infty \b^k<\infty$. Then, the sequence $\{z^k\}$
converges to some optimal solution $\tilde{z}^*\in Z^*$.
\end{Lemma 1}

\section*{Appendix B}
In this section, we establish the relations in Step I, Step II, and Step III of Theorem 3's proof.

 Step I: Relationship for $\frac{2}{L}\left(F(\bar x^{k+1}) -
F(\theta^*)\right)$.

Since $F$ is convex with a Lipschitz gradient, we have:
\[F(y)\le F(x)+\la \nabla F(x),y-x\ra+\frac{L}{2}\|y-x\|^2, \: \forall y,x\in\mathbb{R}^d\]
Letting $y=\bar x^{k+1}$ and $x=\bar x^k$ in the preceding relation
and using  (\ref{eq:bar_y_k}) as well as
$F(\bar{x}^k)=\frac{1}{m}\sum_{i=1}^mf_i(\bar x^k)$, we obtain
\[F(\bar x^{k+1})
\le F(\bar x^k)-\frac{1}{m}\sum_{i=1}^m \la \nabla f_i(\bar
x^k),\bar y^k\ra +\frac{L}{2}\|\bar y^k\|^2\]

Subtracting $F(\theta^*)$ on both sides and multiplying  $2\bar
\lambda^k$ yield
\begin{equation}\label{eq-rel7}
\begin{aligned}
&2\bar\lambda^k \left(F(\bar x^{k+1}) - F(\theta^*)\right) \le
2\bar\lambda^k\left(F(\bar
x^k) - F(\theta^*)\right) \\
&\qquad\qquad-\frac{2\bar\lambda^k}{m}\sum_{i=1}^m \la \nabla
f_i(\bar x^k),\bar y^k\ra +\bar \lambda^k L\|\bar y^k\|^2
\end{aligned}
\end{equation}

The term   $-\frac{2\bar\lambda^k}{m}\sum_{i=1}^m \la \nabla
f_i(\bar x^k),\bar y^k\ra$ satisfies
\begin{equation}\label{eq-estiny}
\begin{aligned}
&-2\left\la \frac{\bar\lambda^k}{m}\sum_{i=1}^m \la \nabla f_i(\bar
x^k),\bar y^k\right\ra\\
&=\left\|\frac{\bar\lambda^k}{m}\sum_{i=1}^m \nabla f_i(\bar x^k) -
\bar y^k\right\|^2
 -\left\|\frac{\bar\lambda^k}{m}\sum_{i=1}^m \nabla f_i(\bar x^k)
\right\|^2 -\left\|\bar y^k\right\|^2
\end{aligned}
\end{equation}
For the first term on the right hand side of  \eqref{eq-estiny},
 by adding and subtracting
$\frac{1}{m}\sum_{i=1}^m \lambda_i^k \nabla f_i(\bar x^k) $, we
obtain
\[
\begin{aligned}
&\left\|\frac{\bar\lambda^k}{m}\sum_{i=1}^m \nabla f_i(\bar x^k) -
\bar y^k\right\|^2 \\
& = \left\|\frac{1}{m}\sum_{i=1}^m (\bar\lambda^k-\lambda_i^k)
\nabla f_i(\bar x^k) +\frac{1}{m}\sum_{i=1}^m \lambda_i^k\nabla
f_i(\bar x^k) - \bar y^k\right\|^2\cr &\le
2\left\|\frac{1}{m}\sum_{i=1}^m (\bar\lambda^k-\lambda_i^k) \nabla
f_i(\bar x^k) \right\|^2 \\
&\qquad +2\left\|\frac{1}{m}\sum_{i=1}^m \lambda_i^k\left(\nabla
f_i(\bar x^k) -\nabla f_i(x_i^k)\right)\right\|^2
\end{aligned}
\]
where we used $\bar y^k=\frac{1}{m}\sum_{i=1}^m \lambda_i^k\nabla
f_i(x_i^k)$ in (\ref{eq:eq-bary}). Using the assumption that each
$\nabla f_i(\cdot)$ is Lipschitz continuous with a constant $L$, we
can further rewrite the preceding inequality as
\begin{equation}\label{eq:StepI_(1)}
\begin{aligned}
&\left\|\frac{\bar\lambda^k}{m}\sum_{i=1}^m \nabla f_i(\bar x^k) -
\bar y^k\right\|^2  \le  \frac{2}{m}\sum_{i=1}^m
(\bar\lambda^k-\lambda_i^k)^2 \|\nabla
f_i(\bar x^k)\|^2 \\
&\qquad +\frac{2}{m}\sum_{i=1}^m (\lambda_i^k)^2 \|\nabla f_i(\bar
x^k) -\nabla f_i(x_i^k)\|^2\cr   &\le
\frac{2}{m}\|\lambda^k-\bar\lambda^k{\bf 1}\|^2\sum_{i=1}^m\|\nabla
f_i(\bar x^k)\|^2 \\
&\qquad+\frac{2L^2}{m}(\lambda_{\max}^k)^2\sum_{i=1}^m \|\bar x^k
-x_i^k\|^2
\end{aligned}
\end{equation}

We next proceed to analyze $\sum_{i=1}^m \|\nabla f_i(\bar
x^k)\|^2$.

By Assumption \ref{assumption:L_and_G}, each $\nabla f_i(\cdot)$ is
Lipschitz continuous with a constant $L$, so we have
\[f_i(v)+\la\nabla f_i(v),u-v\ra+\frac{1}{2L}\|\nabla f_i(v)-\nabla f_i(u)\|^2\le
f_i(u), \forall u,v\] Letting $v=\theta^*$ and $u=\bar x^k$, and
summing the resulting relations over $i=1,\ldots,m$ yield $
F(\theta^*) +  \la  \nabla F(\theta^*), \bar x^k - \theta^*\ra
+\frac{ \sum_{i=1}^m\|\nabla f_i(\theta^*)-\nabla f_i(\bar
x^k)\|^2}{2mL}\le  F (\bar x^k)$. Using $\nabla F(\theta^*)=0$,  we
have
\[\sum_{i=1}^m\|\nabla f_i(\theta^*)-\nabla f_i(\bar x^k)\|^2\le 2mL(F (\bar x^k)-F(\theta^*))\]
Thus, it follows
\begin{equation}\label{eq:nabla_f_i}
\begin{aligned}
&\sum_{i=1}^m\|\nabla f_i(\bar x^k)\|^2\\
&\quad \le  \sum_{i=1}^m 2 \left(\|\nabla f_i(\bar x^k)-\nabla
f_i(\theta^*)\|^2+ \|\nabla f_i(\theta^*)\|^2\right)\cr & \quad \le
4mL(F (\bar x^k) - F(\theta^*))+ 2\sum_{i=1}^m\|\nabla
f_i(\theta^*)\|^2
\end{aligned}
\end{equation}

Combining (\ref{eq-rel7}), (\ref{eq-estiny}), (\ref{eq:StepI_(1)}),
and (\ref{eq:nabla_f_i}) leads to
\[
\begin{aligned}
& 2\bar\lambda^k\left(F(\bar x^{k+1}) - F(\theta^*)\right) \le
2\bar \lambda^k\left(F(\bar x^k) -
F(\theta^*)\right) \\
&\quad+ 8L\|\lambda^k-\bar\lambda^k{\bf 1}\|^2(F (\bar x^k)-F(\theta^*)) \\
&\quad+ \frac{4}{m}\|\lambda^k-\bar\lambda^k{\bf
1}\|^2\sum_{i=1}^m\|\nabla
f_i(\theta^*)\|^2-(\bar\lambda^k)^2\|\nabla F(\bar x^k)\|^2\cr &
\quad+\frac{2L^2}{m}(\lambda_{\max}^k)^2\sum_{i=1}^m \|\bar x^k
-x_i^k\|^2 +(\bar \lambda^k L-1)\|\bar y^k\|^2
\end{aligned}
\] i.e.,
\begin{equation}\label{eq:StepI_final}
\begin{aligned}
& \frac{2}{L}\left(F(\bar x^{k+1}) - F(\theta^*)\right) \le \frac{2}{L}\left(F(\bar x^k) - F(\theta^*)\right)\\
&   + \frac{c_1}{\bar \lambda^k}\|\lambda^k-\bar\lambda^k{\bf
1}\|^2\left(\frac{2}{L}(F (\bar x^k)-F(\theta^*))+
\sum_{i=1}^m\|\nabla
f_i(\theta^*)\|^2\right)\\
&  +\frac{2L}{m\bar\lambda^k}(\lambda_{\max}^k)^2\sum_{i=1}^m
\|\bar x^k -x_i^k\|^2 \\
&  -\frac{\bar\lambda^k}{L}\|\nabla F(\bar x^k)\|^2+\frac{\bar
\lambda^k L-1}{\bar \lambda^k L} \|\bar y^k\|^2
\end{aligned}
\end{equation}
where $c_1=\max\{4L,\, 4/(mL)\}$.

Step II: Relationship for $\sum_{i=1}^m\|x_i^{k+1}-\bar x^{k+1}\|$ and  $\sum_{i=1}^m\|x_i^{k+1} - x_i^k \|^2$.

For the convenience of analysis, we write the iterates of algorithm
PDG-NDS on per-coordinate expressions. Define for all
$\ell=1,\ldots,d,$ and $k\ge0$,
\[
\begin{aligned}
x^k(\ell)&=([x_1^k]_\ell,\ldots,[x_m^k]_\ell)^T,\:
y^k(\ell)=([y_1^k]_\ell,\ldots,[y_m^k]_\ell)^T,\\
g^k(\ell)&=([g_1^k]_\ell,\ldots,[g_m^k]_\ell)^T. \end{aligned}
\] In
this per-coordinate view, (\ref{eq:y_k_definition}) and
(\ref{eq:y_k+1}) has the following form  for all $\ell=1,\ldots,d,$
and $k\ge0$,
\begin{align}\label{eq-alg2-percord}
x^{k+1}(\ell)&=Wx^k(\ell)-y^k(\ell)\cr y^{k+1}(\ell)&=Wy^k(\ell) +
B^k\left(\Lambda^{k+1} g^{k+1}(\ell) -\Lambda^{k}g^k(\ell)\right)
\end{align}

From the definition of $x^{k+1}(\ell)$ in~\eqref{eq-alg2-percord},
and the relation for the average $\bar x^{k+1}$ in
(\ref{eq:bar_y_k}), we obtain for all $\ell=1,\ldots,d$,
\[ x^{k+1}(\ell) - [\bar x^{k+1}]_\ell{\bf 1}
=W\left(x^k(\ell) -[\bar x^k]_\ell{\bf 1}\right) -(y^k(\ell) -[\bar
y^k]_\ell{\bf 1})\] where we use $W{\bf 1}={\bf 1}$. Noting that
$[\bar x^k]_\ell$ is the average of $x^k(\ell)$, i.e.,
$\frac{1}{m}{\bf 1}{\bf 1}^T\left(x^k(\ell) -[\bar x^k]_\ell{\bf
1}\right)=0$, we have
\[
\begin{aligned}
x^{k+1}(\ell) - [\bar x^{k+1}]_\ell{\bf 1} = &\bar W\left(x^k(\ell)
-[\bar x^k]_\ell{\bf 1}\right) -(y^k(\ell) -[\bar y^k]_\ell{\bf 1})
\end{aligned}
\]
where $\bar{W}=W-\frac{ {\bf 1}{\bf 1}^T}{m} $.
 So
it follows
\[\|x^{k+1}(\ell)-[\bar x^{k+1}]_\ell{\bf 1}\|
\le \eta\|x^k-\bar x^k{\bf 1}\| +\|y^k(\ell)-[\bar y^k]_{\ell}{\bf
1}\|\] with $\eta=\|W-\frac{1}{m}{\bf 1}{\bf 1}^T\|<1$. Taking
squares on both sides of the preceding relation, and  using the
inequality $(a+b)^2\le (1+\epsilon)a^2+(1+\epsilon^{-1})b^2$, valid
for any scalars $a$ and $b$, and $\epsilon>0$, we obtain
\[
\begin{aligned}
\|x^{k+1}(\ell)-[\bar x^{k+1}]_\ell{\bf 1}\|^2 \le&
\eta^2(1+\epsilon)\|x^k(\ell)-[\bar x^k]_\ell{\bf 1}\|^2\\
& +(1+\epsilon^{-1})\|y^k(\ell)-[\bar y^k]_\ell{\bf 1}\|^2
\end{aligned}
\] By using  $\eta\in(0,1)$ and
letting $\epsilon=\frac{1-\eta}{\eta}$ which implies
$1+\epsilon=\eta^{-1}$ and $1+\epsilon^{-1}=(1-\eta)^{-1}$, we have
\[
\begin{aligned}
\|x^{k+1}(\ell)-[\bar x^{k+1}]_\ell{\bf 1}\|^2 \le
&\eta\|x^k(\ell)-[\bar x^k]_\ell{\bf 1}\|^2
\\
&+(1-\eta)^{-1}\|y^k(\ell)-[\bar y^k]_\ell{\bf 1}\|^2
\end{aligned}
\] Summing the
preceding relations over $\ell=1,\ldots,d$, and noting
$\sum_{\ell=1}^d\|x^{k+1}(\ell)-[\bar x^{k+1}]_\ell{\bf 1}\|^2
=\sum_{i=1}^m \|x_i^{k+1} -\bar x^{k+1}\|^2$,
$\sum_{\ell=1}^d\|x^k(\ell)-[\bar x^k]_\ell{\bf
1}\|^2=\sum_{i=1}^m\|x_i^k-\bar x^k\|^2$, and
$\sum_{\ell=1}^d\|y^k(\ell)-[\bar y^k]_\ell{\bf
1}\|^2=\sum_{i=1}^m\|y_i^k-\bar y^k\|^2$, we obtain
\begin{equation}\label{eq:StepII_final_a}
\begin{aligned}
 \sum_{i=1}^m \|x_i^{k+1} -\bar x^{k+1}\|^2 \le& \eta\sum_{i=1}^m
\|x_i^{k+1} -\bar x^{k+1}\|^2
\\
&+(1-\eta)^{-1}\sum_{i=1}^m\|y_i^k-\bar y^k\|^2
\end{aligned}
\end{equation}

Next we proceed to analyze $\sum_{i=1}^m\|x_i^{k+1} - x_i^k \|^2$.
Using (\ref{eq:y_k_definition}), we have for every coordinate index
$\ell=1,\ldots,d$
\[
\begin{aligned}
&x^{k+1}(\ell)-x^k(\ell) =Wx^k(\ell)-y^k(\ell)
-x^k(\ell)\\
&=(W-I)x^k(\ell)-y^k(\ell) =(W-I)(x^k(\ell)-[\bar x^k]_\ell{\bf
1})-y^k(\ell)
\end{aligned}
\] where we used the fact
$(W-I){\bf 1}=0$. By letting $r=\|W-I\|$, we obtain
\[
\begin{aligned}
&\|x^{k+1}(\ell)-x^k(\ell)\| \le r \|x^k(\ell)-[\bar x^k]_\ell{\bf
1}\|+\|y^k(\ell)\| \\
&\qquad \le r\|x^k(\ell)-[\bar x^k]_\ell{\bf 1}\|+\|y^k -[ \bar
y^k]_\ell{\bf 1}\| +\sqrt{m}|[\bar y^k]_\ell |
\end{aligned}
\] where the last inequality is obtained by adding and
subtracting $[\bar y^k]_\ell{\bf 1}$ to $y^k(\ell)$, and using the
triangle inequality for the norm. Thus, we have {
\[
\begin{aligned}
\|x^{k+1}(\ell)-x^k(\ell)\|^2 \le& 3r^2\|x^k(\ell)-[\bar
x^k]_\ell{\bf 1}\|^2+\\
& 3\|y^k -[ \bar y^k]_\ell{\bf 1}\|^2 + 3m|[\bar y^k]_\ell |^2
\end{aligned}
\] }

By summing over  $\ell=1,\ldots,d$, we obtain
\begin{equation}\label{eq:StepII_final_b}
 \begin{aligned}
 \sum_{i=1}^m\|x_i^{k+1} - x_i^k \|^2 \le& 3r^2\sum_{i=1}^m\|x_i^k -
\bar x^k \|^2 \\
&+3\sum_{i=1}^m\|y_i^k - \bar y^k\|^2 + 3m\|\bar y^k\|^2
 \end{aligned}
\end{equation}

Step III: Relationship for $\sum_{i=1}^m\|y_i^{k+1}-\bar y^{k+1}{\bf
1}\|$.

Using the column stochastic property of $B^k$, from
(\ref{eq:y_k+1}), the $\ell$th entries of  $[\bar y^k]_\ell$
satisfy
\[[\bar y^{k+1}]_\ell
=[\bar y^k]_\ell + \frac{1}{m}{\bf 1}^T\Lambda^{k+1} g^{k+1}(\ell)-
\frac{1}{m}{\bf 1}^T\Lambda^k g^k(\ell)\] Then, using
 \eqref{eq-alg2-percord}, we obtain for all $\ell=1,\ldots,d$,
\[
\begin{aligned}
y^{k+1}(\ell) - [\bar y^{k+1}]_\ell{\bf 1}  =& \bar W(y^k(\ell)-[\bar y^k]_\ell{\bf 1}) +\\
&  \bar B^k\Lambda^{k+1}g^{k+1}(\ell) -\bar B^k\Lambda^k g^k(\ell)
\end{aligned}
\]
where $\bar W=W-\frac{1}{m} {\bf 1}{\bf 1}^T$ and $\bar
B^k=\left(B^k-\frac{1}{m}{\bf 1}{\bf 1}^T\right)$.

By adding and subtracting $ \bar B^k\Lambda^{k+1}g^k(\ell)$, and
taking the Euclidean norm, we find that for all $\ell=1,\ldots,d$,
\begin{equation}\label{eq:stepIII_1}
\begin{aligned}
& \|y^{k+1}(\ell) -[\bar y^{k+1}]_\ell{\bf 1}\|\le \eta
\|y^k(\ell) - [\bar y^k]_\ell{\bf 1}\|\\
& + \tau \left\|\Lambda^{k+1}\left(g^{k+1}(\ell) -
g^k(\ell)\right)\right\| + \tau
\left\|\left(\Lambda^{k+1}-\Lambda^k\right)g^k(\ell)\right\|
\end{aligned}
\end{equation}
where $\eta=\|\bar W\|$ and $\tau=\|\bar B^k\|$.

Since we always have
\[\left\|\bar B^k\right\|\le \left\|\bar B^k\right\|_F
 \leq m
\]
where $\|\cdot\|_F$ denotes the Frobenius matrix norm, we have
\begin{equation}\label{eq:tau_m}
\tau\le   m
\end{equation}

Using the fact that
$\Lambda^k$ is a diagonal matrix for all $k\geq 0$, i.e., $\Lambda^k={\rm
diag}(\lambda^k)$, we have
\begin{equation}\label{eq:stepIII_3}
\begin{aligned}
&\left\|\Lambda^{k+1}\left(g^{k+1}(\ell)-g^k(\ell)\right)\right\|
\le
\|\lambda^{k+1}\|\,\|g^{k+1}(\ell)-g^k(\ell)\|,\\
&\left\|\left(\Lambda^{k+1}-\Lambda^k\right)g^k(\ell)\right\|\le
\|\lambda^{k+1}-\lambda^k\|\,\|g^k(\ell)\|
\end{aligned}
\end{equation}

Therefore,
combining (\ref{eq:stepIII_1}), (\ref{eq:tau_m}), and (\ref{eq:stepIII_3}) leads to
\[
\begin{aligned}
&\|y^{k+1}(\ell) - [\bar y^{k+1}]_\ell{\bf 1}\|  \le  \eta
\|y^k(\ell) - [\bar y^k]_\ell{\bf 1}\|+\\
&m \|\lambda^{k+1}\|\,\|g^{k+1}(\ell)-g^k(\ell)\| + m
\left\|\lambda^{k+1}-\lambda^k\right\|\,\|g^k(\ell)\|
\end{aligned}
\]
Thus, by taking squares and using $(a+b)^2\le(1+\epsilon)a^2
+(1+\epsilon^{-1})b^2$ holding for any $\epsilon>0$, we obtain the
following inequality by setting $\epsilon=\frac{1-\eta}{\eta}$:
\[
\begin{aligned}
&\|y^{k+1}(\ell) - [\bar y^{k+1}]_\ell{\bf 1}\|^2 \le \eta
\|y^k(\ell) - [\bar y^k]_\ell{\bf 1}\|^2+\frac{2m^2}{1-\eta}\times\\
&\left(\|\lambda^{k+1}\|^2\,\|g^{k+1}(\ell)-g^k(\ell)\|^2
+\left\|\lambda^{k+1}-\lambda^k\right\|^2\,\|g^k(\ell)\|^2\right)
\end{aligned}
\]
By summing these relations over $\ell=1,\ldots,d$, we find
\begin{equation}\label{eq-y2}
\begin{aligned}
& \sum_{i=1}^m \|y^{k+1}_i - \bar y^{k+1}\|^2 \le \eta
\sum_{i=1}^m\|y^k_i - \bar y^k\|^2+\frac{2m^2}{1-\eta}\times
\\
&\left(\|\lambda^{k+1}\|^2\,\sum_{i=1}^m\|g^{k+1}_i-g^k_i\|^2
+\left\|\lambda^{k+1}-\lambda^k\right\|^2\,\sum_{i=1}^m\|g^k_i\|^2\right)
\end{aligned}
\end{equation}

Next we bound $\sum_{i=1}^m\|g^{k+1}_i-g^k_i\|^2$ and
$\sum_{i=1}^m\|g^k_i\|^2$.

Since every $\nabla f_i$ is Lipschitz continuous with $L>0$, we have
\[\sum_{i=1}^m \|g^{k+1}_i - g^k_i\|^2\le L^2 \sum_{i=1}^m \|x^{k+1}_i-x^k_i\|^2\]
which, in combination with (\ref{eq:StepII_final_b}), leads to
\begin{equation}\label{eq:stepIII_4}
\begin{aligned}
\sum_{i=1}^m \|g^{k+1}_i - &g^k_i\|^2 \le 3r^2L^2 \sum_{i=1}^m
\|x^k_i-\bar x^k\|^2 \\
&+3L^2\sum_{i=1}^m\|y^k_i- \bar y^k\|^2+3mL^2\|\bar y^k\|^2
\end{aligned}
\end{equation}
Using the Lipschitz continuity of each $g^k_i$, we obtain
\[
\begin{aligned}
\sum_{i=1}^m\|g_i^k\|^2&=\sum_{i=1}^m\|\nabla f_i(x_i^k)\|^2
\\&\le 2L^2
\sum_{i=1}^m\|x_i^k - \bar x^k\|^2 +2 \sum_{i=1}^m \|\nabla f_i(\bar
x^k)\|^2
\end{aligned}
\]
which, in combination with (\ref{eq:nabla_f_i}), leads to
\begin{equation}\label{eq:StepIII_5}
\begin{aligned}
\sum_{i=1}^m\|g_i^k\|^2 \le& 2L^2 \sum_{i=1}^m\|x_i^k - \bar x^k\|^2
+8mL(F (\bar x^k)-F(\theta^*))\\
& +4\sum_{i=1}^m\|\nabla f_i(\theta^*)\|^2
\end{aligned}
\end{equation}
By substituting (\ref{eq:stepIII_4}) and (\ref{eq:StepIII_5}) into
(\ref{eq-y2}), we   obtain
\begin{equation}\label{eq:stepIII_final}
\begin{aligned}
\hspace{-0.1cm} &\sum_{i=1}^m \|y_i^{k+1}-\bar y^{k+1}\|^2\le  (\eta
+\frac{6m^2L^2}{1-\eta}\|\lambda^{k+1}\|^2 )\sum_{i=1}^m\|y_i^k-\bar
y^k\|^2\cr \hspace{-0.1cm}&
+\frac{6m^2r^2L^2}{1-\eta}\left(\left\|\lambda^{k+1}\right\|^2\sum_{i=1}^m\|x_i^k-\bar
x^k\|^2 + m\left\|\lambda^{k+1}\right\|^2\|\bar y^k\|^2\right)\cr
\hspace{-0.1cm}& +
\frac{2m^2}{1-\eta}\left\|\lambda^{k+1}-\lambda^k\right\|^2\,
\left(2L^2 \sum_{i=1}^m\|x_i^k - \bar x^k\|^2 +\right.\\
\hspace{-0.1cm}&\qquad\qquad\qquad \left.8mL(F (\bar
x^k)-F(\theta^*))+ 4\sum_{i=1}^m\|\nabla f_i(\theta^*)\|^2\right)
\end{aligned}
\end{equation}

\bibliographystyle{unsrt}

\bibliography{reference1}
\vspace{-1.2cm}

\end{document}